\documentclass[11pt,letterpaper]{article}

\usepackage[title]{appendix}
\usepackage[english]{babel}
\usepackage[utf8]{inputenc}
\usepackage{amsmath}
\usepackage{mathtools}
\usepackage{amssymb}
\usepackage{upgreek}
\usepackage{graphicx}
\usepackage{geometry}
\usepackage{cite}
\usepackage{authblk}
\usepackage{color}
\usepackage{soul}
\usepackage{xcolor}

\usepackage{float}
\usepackage{multirow}
\usepackage{booktabs}
\usepackage{tabu} 
\usepackage{xcolor}
\usepackage{dcolumn}
\newcolumntype{A}{D{.}{.}{2.3}}
\usepackage{hyperref}
\title{Design Optimization of Dynamic Flexible Multibody Systems Using the Discrete Adjoint Variable Method\thanks{This is a preprint of an article published in \emph{Computers \& Structures}. The final authenticated version is available online at: \href{https://doi.org/10.1016/j.compstruc.2018.12.007}{https://doi.org/10.1016/j.compstruc.2018.12.007}.}}

\author[a]{\small Mehran Ebrahimi\footnote{\textit{E-mail address: }mehran.ebrahimi@autodesk.com}\textsuperscript{,}}
\author[a]{\small Adrian Butscher}
\author[a]{\small Hyunmin Cheong}
\author[a]{\small Francesco Iorio}

\affil[a]{\footnotesize Autodesk Research, 661 University Avenue, Toronto, ON M5G 1M1, Canada}

\date{}

\begin{document}
\maketitle

\begin{abstract}
\noindent
The design space of dynamic multibody systems (MBSs), particularly those with flexible components, is considerably large. Consequently, having a means to efficiently explore this space and find the optimum solution within a feasible time-frame is crucial. It is well-known that for problems with several design variables, sensitivity analysis using the adjoint variable method extensively reduces the computational costs. This paper presents the novel extension of the discrete adjoint variable method to the design optimization of dynamic flexible MBSs. The extension involves deriving the adjoint equations directly from the discrete, rather than the continuous, equations of motion. This results in a system of algebraic equations that is computationally less demanding to solve compared to the system of differential algebraic equations produced by the continuous adjoint variable method. To describe the proposed method, it is integrated with a numerical time-stepping algorithm based on geometric variational integrators. The developed technique is then applied to the optimization of MBSs composed of springs, dampers, beams and rigid bodies, considering both geometrical (e.g., positions of joints) and non-geometrical (e.g., mechanical properties of components) design variables. To validate the developed methods and show their applicability, three numerical examples are provided.
\end{abstract}

\textit{Keywords} Flexible multibody systems, Design optimization, Sensitivity analysis, Adjoint variable method, Discrete adjoint variable method, Geometrical design variables

\section{Introduction}
\label{sec:introduction}
Multibody systems are mechanical assemblies composed of several interconnected parts and have numerous applications in automotive, aerospace, robotics and many other industries. Generally, these systems have a large number of degrees of freedom, show dynamic behaviors, and their components may undergo both large overall motion and large deformation, thus making the optimal design of such systems a challenging task.

The optimization of an MBS involves finding the optimum values of the system's parameters, also called design variables, that minimize (or maximize) a desired objective and satisfy a set of constraints. Design variables, in general, could be classified into control parameters and design parameters. The former are the focus of optimal control problems and their examples include time-dependent input forces and torques. Design parameters, on the other hand, are handled in design optimization problems and are related to mechanical and physical properties of a system and its components (e.g., material properties, cross-sectional areas, lengths). In both optimization scenarios, for an assembly with many coupled flexible and rigid bodies, the parametric design space is considerably vast. Therefore, an efficient algorithm is required to search this space and find an optimum design within a reasonable timescale.  

Accordingly, gradient-based optimizers have been shown to outperform meta-heuristic optimization techniques and have a better convergence rate (see e.g. \cite{Tromme2013}). In order to utilize gradient-based approaches, sensitivity analysis should be performed, which requires computing the derivative of the objective and constraint functions with respect to the design variables. In the majority of problems, these functions depend not only explicitly on the design variables, but also implicitly on them through the state variables, which describe the state of an MBS. Hence, to obtain their total derivatives, one needs to calculate the sensitivity values of the state variables to the design variables as well. This becomes a laborious task, particularly in large-scale and time-dependent problems.

One way of computing the gradients is to calculate them numerically using the finite difference method.  Although this approach is easy to implement, it provides only an approximation to the actual gradient values, and the perturbation step of each design variable is not known \emph{a priori}. Additionally, to obtain the sensitivities, one additional simulation per design variable is required. Consequently, for computationally expensive problems, this method is highly inefficient and prohibitively time-consuming \cite{Greene1989,Maly1996}. A more practical and accurate approach is to apply analytical techniques that compute the exact gradients. 

The two most widely used methods in the analytical category are the \emph{direct differentiation method} and the \emph{adjoint variable method}. In the former, first, the motion equations of a given MBS are differentiated with respect to each design variable. This results in a system of differential-algebraic equations (DAEs), in which the gradients of the state variables and the Lagrange multipliers are unknowns. Then, by solving this system and implementing the solution values into the sensitivity equations of the objective and constraint functions, the desired gradients are evaluated{\mbox{\cite{Maly1996, Mukherjee2008,Dias1997, Wang2005, Haug1982, Pi2012,tromme2015structural}}}. Since for each design variable a new set of DAEs needs to be solved, the computational cost of this approach becomes significant for complex assemblies with a large number of design variables.

An alternative strategy is to exploit the adjoint variable method  \cite{Schaffer2005, Ding2007, Zhu2015, Cao2002, Li2000, Pi2012}. In this technique, the explicit computation of the gradients of state variables and Lagrange multipliers, as required in the direct differentiation method, is avoided by introducing a set of auxiliary variables called \emph{adjoint variables}. To compute the adjoint variables, the equations of motion are first solved (\emph{forward simulation}), followed by integrating a system of linear DAEs known as adjoint equations backward in time (\emph{backward simulation}). The main advantage of the adjoint variable method is that the generated linear system needs to be solved only once at each iteration of the optimization process, while eliminating entirely the need for computing the derivatives of state variables and Lagrange multipliers with respect to the design variables.

In one view, the adjoint variable method can further be split into two categories: the continuous adjoint variable method (CAVM) and the discrete adjoint variable method (DAVM). In the former, the continuous equations of motion (i.e., the equations of motion prior to discretization) are used to derive the continuous adjoint equations --- a system of linear DAEs. This system should be then discretized in time and solved to find the adjoint variables. On the other hand, in the DAVM, the discrete adjoint equations are directly derived from the discretized equations of motion. This approach yields the exact gradients of the discrete objective and constraint functions. More importantly, the equations generated by the DAVM no longer need an additional level of discretization in time and form a system of linear algebraic equations, which are computationally easier to solve compared to the DAEs that arise via the CAVM.

The DAVM has been previously used in optimal control and design optimization problems \cite{ober2011discrete, Lauss2017a}, however, its application in large-scale MBSs with coupled flexible and rigid components is not yet well studied. This paper intorduces the novel application of the DAVM to the design optimization of flexible MBSs. Without loss of generality, to describe the implementation details and solve the equation of motion, the proposed sensitivity analysis scheme is applied to a symplectic-momentum preserving geometric variational integrator proposed by Leyendecker et al. \cite{Leyendecker2008}. The same method can be applied to other numerical time-stepping solvers.

{In the current study, both \emph{geometrical} and \emph{non-geometrical} design variables are taken into account. To the authors' knowledge, the inclusion of geometrical design variables in optimizing dynamic MBSs has not been investigated before. Geometrical design variables concern the parameters defining the shape of a multibody configuration for example the global initial positions of joints and components' length. Non-geometrical variables include other mechanical and physical properties of the bodies such as spring constants, damping coefficients, masses, cross-sections and Young's moduli.  Considering geometrical parameters would increase the dimension of the design space and provides the optimization routine with more flexibility. There are many applications, for instance designing an MBS whose components should follow certain trajectories in space, where incorporating only the non-geometrical parameters would not lead to a solution, while including the geometrical variables is essential to finding the desired design.}

To demonstrate how to apply the proposed methods, without loss of generality, the relevant sensitivity equations for assemblies made of springs, dampers, beams and rigid bodies connected via spherical (pin) and welded (fixed) joints are derived. The presented approach can be extended to include other component and joint types. The remainder of the paper is organized as follows. Section \ref{sec:problemdef} introduces the general problem definition and basics of gradient-based optimization algorithms. Section \ref{sec:Dynamicsof} presents the equations of motion and how to solve them using the selected geometrical variational integrator. In Section \ref{sec:sensitivity}, the detailed derivation of the DAVM is provided. In Section \ref{sec:SensitivityGeo}, first, the equations involved in describing the dynamics of interconnected springs, damper, beams and rigid bodies are discussed. Then, the equations for computing the sensitivities with respect to geometrical design variables are developed. Finally, in the last section, three numerical examples are presented.

\newpage
\section{General problem definition}
\label{sec:problemdef}

Consider an unconstrained optimization problem for the vector $\vec{\mathbf{a}} \in \mathbb{R}^n$ of \emph{design variables} of a multibody system in the following form:

\begin{equation}\label{eq:generalform}
\begin{gathered}
\min_{\vec{\mathbf{a}}} \phi \left(\vec{\mathbf{q}}\left(\vec{\mathbf{a}} \right), \dot{\vec{\mathbf{q}}}\left(\vec{\mathbf{a}} \right), \vec{\boldsymbol{\uplambda}}\left(\vec{\mathbf{a}} \right), \vec{\mathbf{a}} \right) \\
\end{gathered}
\end{equation}
In this equation, $\phi$ is the objective function, $\vec{\mathbf{q}}$ represents the vector of \emph{state variables}, $\dot{\vec{\mathbf{q}}}$ is its \emph{time derivative}, and $\vec{\boldsymbol{\uplambda}}$ denotes the vector of \emph{Lagrange multipliers} associated with the joint (constraint) equations. These quantities are time-dependent and are the solutions of the multibody dynamic equations, to be described in Section \ref{sec:Dynamicsof}. As the behavior of an MBS is affected by the choices of the design variables, $\vec{\mathbf{q}}$, $\dot{\vec{\mathbf{q}}}$ and $\vec{\boldsymbol{\uplambda}}$ are also dependent on $\vec{\mathbf{a}}$. It is assumed that $\phi$ is at least a twice-differentiable function of its arguments.  Moreover, It is assumed that $\phi$ has the following generic form
\begin{equation}
\label{eq:genericform}
\phi = F \left( \vec{\mathbf{q}}_0, \dot{\vec{\mathbf{q}}}_0, \vec{\mathbf{q}}_T, \dot{\vec{\mathbf{q}}}_T, \vec{\mathbf{a}} \right) + \int_0^T H \left( \vec{\mathbf{q}}, \dot{\vec{\mathbf{q}}}, \vec{\boldsymbol{\uplambda}}, \vec{\mathbf{a}} \right) dt
\end{equation}
where $F$ is a function defined on the two ends of the simulation duration, $t=0$ and $t=T$, and $H$ is described over its entire duration.

A gradient-based algorithm for finding a locally minimizing solution of Equation \ref{eq:generalform} proceeds by iterative improvement of the optimization objective function: if $k \geq 0$ denotes the iteration counter, then one strives to improve the current values of the design variables $\vec{\mathbf{a}}^k$ by taking an appropriate small step, i.e., $\vec{\mathbf{a}}^{k+1} := \vec{\mathbf{a}}^k + \delta \vec{\mathbf{a}}^k$.  To determine $\delta \vec{\mathbf{a}}^k$, $\phi$ is represented as the following Taylor series near $\vec{\mathbf{a}}^k$

\begin{equation}\label{eq:taylor}
\begin{aligned}
\phi \left(\vec{\mathbf{q}}, \dot{\vec{\mathbf{q}}}, \vec{\boldsymbol{\uplambda}}, \vec{\mathbf{a}}^k+\delta \vec{\mathbf{a}}^k \right) &= \phi \left(\vec{\mathbf{q}}, \dot{\vec{\mathbf{q}}}, \vec{\boldsymbol{\uplambda}}, \vec{\mathbf{a}}^k \right)+\frac{d}{d \vec{\mathbf{a}}} \phi \left(\vec{\mathbf{q}}, \dot{\vec{\mathbf{q}}}, \vec{\boldsymbol{\uplambda}}, \vec{\mathbf{a}}^k \right)\cdot \delta \vec{\mathbf{a}}^k + \mathrm{o}\left(\left\|\delta \vec{\mathbf{a}}^k \right\|^2 \right) \\
&\mathrm{with} \: \lim_{\delta \vec{\mathbf{a}}^k \rightarrow 0} \frac{\mathrm{o}\left(\left\|\delta \vec{\mathbf{a}}^k \right\|^2 \right)}{\left\|\delta \vec{\mathbf{a}}^k \right\|} = 0
\end{aligned}
\end{equation}

where $d \phi / d \vec{\mathbf{a}}$ is the \emph{sensitivity} of $\phi$ with respect to the design variables. To improve the value of $\phi$ to first order, $\delta \vec{\mathbf{a}}^k$ must be along a \emph{descent direction} such that the second term in Equation \ref{eq:taylor} is negative. One possible choice for such a direction is

\begin{equation}\label{eq:descentdir}
\delta \vec{\mathbf{a}}^k = -\epsilon D_{\vec{\mathbf{a}}} \phi \left(\vec{\mathbf{q}}, \dot{\vec{\mathbf{q}}}, \vec{\boldsymbol{\uplambda}}, \vec{\mathbf{a}}^k \right)
\end{equation}

{in which $D_{\vec{\mathbf{a}}} \phi=\frac{d \phi}{d \vec{\mathbf{a}}}$ and $\epsilon > 0$ is a small descent step.} Putting Equation \ref{eq:descentdir} in Equation \ref{eq:taylor} leads to
\begin{equation}
\phi \left(\vec{\mathbf{q}}, \dot{\vec{\mathbf{q}}}, \vec{\boldsymbol{\uplambda}}, \vec{\mathbf{a}}^k+\delta \vec{\mathbf{a}}^k \right) \approx \phi \left(\vec{\mathbf{q}}, \dot{\vec{\mathbf{q}}}, \vec{\boldsymbol{\uplambda}}, \vec{\mathbf{a}}^k \right)- \epsilon \left\| D_{\vec{\mathbf{a}}} \phi \left(\vec{\mathbf{q}}, \dot{\vec{\mathbf{q}}}, \vec{\boldsymbol{\uplambda}}, \vec{\mathbf{a}}^k \right) \right\|^2
\end{equation}
which guarantees the decrease of $\phi$ at iteration $k$ of the optimization process. This procedure can be easily extended to constrained optimization problems, using for example the Augmented Lagrangian method. Thus, in order to move along the descent direction in Equation \ref{eq:descentdir}, the derivative of $\phi$ with respect to $\vec{\mathbf{a}}$ is desired. It is expressed as
\begin{equation}\label{eq:dphida}
\frac{d \phi}{d \vec{\mathbf{a}}}=\frac{\partial \phi}{\partial \vec{\mathbf{a}}} + \frac{\partial \phi}{\partial \vec{\mathbf{q}}} \frac{\partial \vec{\mathbf{q}}}{\partial \vec{\mathbf{a}}} + \frac{\partial \phi}{\partial \dot{\vec{\mathbf{q}}}} \frac{\partial \dot{\vec{\mathbf{q}}}}{\partial \vec{\mathbf{a}}} + \frac{\partial \phi}{\partial \vec{\boldsymbol{\uplambda}}} \frac{\partial \vec{\boldsymbol{\uplambda}}}{\partial \vec{\mathbf{a}}}
\end{equation}
To compute $\frac{d \phi}{d \vec{\mathbf{a}}}$, the values of $\frac{\partial \vec{\mathbf{q}}}{\partial \vec{\mathbf{a}}}$, $\frac{\partial \dot{\vec{\mathbf{q}}}}{\partial \vec{\mathbf{a}}}$ and $\frac{\partial \vec{\boldsymbol{\uplambda}}}{\partial \vec{\mathbf{a}}}$, as well as $\frac{\partial \phi}{\partial \vec{\mathbf{a}}}$, for different types of design variables are required. This is the main focus of the present paper and is addressed in the forthcoming sections. But first, the following section presents the equations of motions for an MBS and the method chosen to solve them.

\section{Solving the equations of motion}
\label{sec:Dynamicsof}

According to \emph{Hamilton's principle} (a.k.a.\ the \emph{least action principle}), it is possible to re-formulate Newton's Law of Motion in such a way that the trajectories of a dynamic system themselves satisfy a type of optimization problem.  Namely, the trajectory that an MBS takes to move between two positions in space minimizes a quantity known as the \emph{action integral}. Therefore, the motion of an MBS can be characterized by finding the stationary solutions of its action integral. This leads to a useful method of discretizing the equations of motion known as \emph{geometric variational integration}, which is outlined below. 

Suppose $\vec{\mathbf{q}} \in \mathbb{R}^m $ is the vector of all degrees of freedom of an MBS and $\dot{\vec{\mathbf{q}}}$ is their time derivative. Also, let $\vec{\boldsymbol{\uplambda}} \in \mathbb{R}^l $ denotes the vector of Lagrange multipliers associated with $l$ \emph{holonomic} constraints (joints) between different bodies in the assembly. Assuming a \emph{conservative} system, its continuous action integral for time $t \in \left[0,T \right] $ is expressed by
\begin{equation}\label{eq:action}
S (\vec{\mathbf{q}}) = \int_0^T \left( L \left( \vec{\mathbf{q}}, \dot{\vec{\mathbf{q}}} \right) -
\vec{\boldsymbol{\uplambda}}\cdot \vec{\mathbf{g}} (\vec{\mathbf{q}})
 \right) dt
\end{equation}
where $L \left( \vec{\mathbf{q}}, \dot{\vec{\mathbf{q}}} \right)$ and $\vec{\mathbf{g}} (\vec{\mathbf{q}})$ are the Lagrangian and constraint equations, respectively. The Lagrangian function $L$ is defined as the \emph{kinetic energy} minus the \emph{potential energy} of the system
\begin{equation}\label{eq:lagrangian}
L \left( \vec{\mathbf{q}}, \dot{\vec{\mathbf{q}}} \right) = T(\dot{\vec{\mathbf{q}}})-U(\vec{\mathbf{q}}) =  \frac{1}{2} \dot{\vec{\mathbf{q}}}^{\mathrm{T}} \mathbf{M} \dot{\vec{\mathbf{q}}} - U(\vec{\mathbf{q}})
\end{equation}
in which $\mathbf{M}$ is the positive-definite \emph{mass matrix} of the system, and $U$ is the total potential energy, due to for example gravity and the deformation of springs and flexible bodies. Following Hamilton's principle, two classes of time-stepping algorithms have come to exist: 1) computing the stationary values of the continuous action integral and then discretizing the resultant dynamic equations in time, or 2) first discretizing the continuous action integral in time and then deriving the discrete dynamic equations directly from the discretized action integral \cite{Marsden2001}. The latter approach forms the basis of geometric variational integration and is pursued in this paper.

\subsection{Geometric variational integrators}
\label{subsec:geometric}

Geometric variational integrators belong to the class of \emph{symplectic-momentum} preserving techniques. They are applicable to a wide spectrum of dynamic problems and have been proven superior to traditional approaches \cite{Kane2000, Stern2006}. Traditional time-stepping algorithms (such as explicit Euler, implicit Euler, Runge-Kutta) suffer from numerical instabilities and artificial dissipation, which makes them incapable of capturing the true dynamic behavior of a system, particularly in the long-duration problems \cite{Kane2000}. This can have significant consequences on the accuracy of \emph{forward} (solving the dynamic equations) and \emph{backward} (solving the adjoint equations) simulations for the system of interest. In the current work, a specific type of variational integrators proposed by Leyendecker et al.\ \cite{Leyendecker2008} is adopted to solve the equations of motion and further describe the proposed sensitivity computation scheme.

Following the idea of variational integrators, if the time domain $[0,T]$ is split into $N$ intervals as $[t_n,t_n+h_n ]  (n=0,1,\dots,N-1)$, Equation \ref{eq:action} can be rewritten as
\begin{equation}\label{eq:actiondis}
\int_0^T \left( L \left( \vec{\mathbf{q}}, \dot{\vec{\mathbf{q}}} \right) -
\vec{\boldsymbol{\uplambda}}\cdot \vec{\mathbf{g}} (\vec{\mathbf{q}})
 \right) dt = \displaystyle\sum_{n=0}^{N-1} \left(\int_{t_n}^{t_{n+1}} \left( L \left( \vec{\mathbf{q}}, \dot{\vec{\mathbf{q}}} \right) -
\vec{\boldsymbol{\uplambda}}\cdot \vec{\mathbf{g}} (\vec{\mathbf{q}})
 \right) dt \right) 
\end{equation}
The following quadrature approximations can be applied to the integral summands:
\begin{equation}\label{eq:quadrature}
\begin{aligned}
\int_{t_n}^{t_{n+1}} \left( L \left( \vec{\mathbf{q}}, \dot{\vec{\mathbf{q}}} \right) \right) dt &\approx 
h_n  L \left( (1-\alpha)\vec{\mathbf{q}}_n + \alpha \vec{\mathbf{q}}_{n+1}, \frac{\vec{\mathbf{q}}_{n+1}-\vec{\mathbf{q}}_n}{h_n} \right) = L_d \left( \vec{\mathbf{q}}_{n},\vec{\mathbf{q}}_{n+1} \right) \\
\int_{t_n}^{t_{n+1}} \left( \vec{\boldsymbol{\uplambda}}\cdot\vec{\mathbf{g}}(\vec{\mathbf{q}}) \right) dt &\approx 
\frac{h_n}{2} \left( \vec{\boldsymbol{\uplambda}}_n\cdot\vec{\mathbf{g}}(\vec{\mathbf{q}}_n) + \vec{\boldsymbol{\uplambda}}_{n+1}\cdot\vec{\mathbf{g}}(\vec{\mathbf{q}}_{n+1}) \right) = \frac{1}{2} \left( \vec{\boldsymbol{\uplambda}}_n\cdot\vec{\mathbf{g}}_d(\vec{\mathbf{q}}_n) + \vec{\boldsymbol{\uplambda}}_{n+1}\cdot\vec{\mathbf{g}}_d (\vec{\mathbf{q}}_{n+1}) \right)
\end{aligned}
\end{equation}
where $\alpha \in [0,1]$. If $\alpha=0.5$, the approximations are second-order accurate, otherwise they are of linear accuracy. Higher order of accuracies can be achieved by improving the quadrature rule. Using Equations \ref{eq:actiondis} and \ref{eq:quadrature}, the discretized action integral is given by
\begin{equation}\label{eq:actiondis2}
S_d(\vec{\mathbf{q}}) = \displaystyle\sum_{n=0}^{N-1} \left( L_d \left( \vec{\mathbf{q}}_n,\vec{\mathbf{q}}_{n+1} \right) - \frac{1}{2} \left( \vec{\boldsymbol{\uplambda}}_n\cdot\vec{\mathbf{g}}_d(\vec{\mathbf{q}}_n) + \vec{\boldsymbol{\uplambda}}_{n+1}\cdot\vec{\mathbf{g}}_d (\vec{\mathbf{q}}_{n+1}) \right) \right)
\end{equation}
Taking the variation of this equation, re-indexing it, and setting $\delta \vec{\mathbf{q}}_0 = \delta \vec{\mathbf{q}}_N=0$, lead to
\begin{equation}\label{eq:deltaactiondis2}
\delta S_d(\vec{\mathbf{q}}) = \displaystyle\sum_{n=1}^{N-1} \left( \frac{\partial L_d \left( \vec{\mathbf{q}}_n,\vec{\mathbf{q}}_{n+1} \right)}{\partial\vec{\mathbf{q}}_n} + \frac{\partial L_d \left( \vec{\mathbf{q}}_{n-1},\vec{\mathbf{q}}_n \right)}{\partial\vec{\mathbf{q}}_n} - \left[\frac{\partial \vec{\mathbf{g}}_d (\vec{\mathbf{q}}_n)}{\partial \vec{\mathbf{q}}_n} \right]^\mathrm{T} \vec{\boldsymbol{\uplambda}}_n
\right)\cdot \delta \vec{\mathbf{q}}_n
\end{equation}
The stationary values of Equation \ref{eq:deltaactiondis2}, along with the constraint equations, provide the \emph{discrete Euler-Lagrange} equation
\begin{equation}\label{eq:eulerLagrange}
\begin{aligned}
\frac{\partial L_d \left( \vec{\mathbf{q}}_n,\vec{\mathbf{q}}_{n+1} \right)}{\partial\vec{\mathbf{q}}_n} + \frac{\partial L_d \left( \vec{\mathbf{q}}_{n-1},\vec{\mathbf{q}}_n \right)}{\partial\vec{\mathbf{q}}_n} - \left[\frac{\partial \vec{\mathbf{g}}_d (\vec{\mathbf{q}}_n)}{\partial \vec{\mathbf{q}}_n} \right]^\mathrm{T} \vec{\boldsymbol{\uplambda}}_n &= \vec{\mathbf{0}} \\
\vec{\mathbf{g}}(\vec{\mathbf{q}}_{n+1}) &= \vec{\mathbf{0}}
\end{aligned}
\end{equation}
Equation \ref{eq:eulerLagrange} is a system of nonlinear algebraic equations that can be solved iteratively, using for example the Newton-Raphson method, to find $\vec{\mathbf{q}}_{n+1}$ and $\vec{\boldsymbol{\uplambda}}_n$ for $n=1,2,\dots,N-1$. Note: in the presence of \emph{non-conservative} external forces, Equation \ref{eq:eulerLagrange} is modified to the \emph{discrete Lagrange-d'Alembert} equation \cite{Kane2000}

\begin{equation}\label{eq:lagrangeDalembert}
\begin{aligned}
\frac{\partial L_d \left( \vec{\mathbf{q}}_n,\vec{\mathbf{q}}_{n+1} \right)}{\partial\vec{\mathbf{q}}_n} + \frac{\partial L_d \left( \vec{\mathbf{q}}_{n-1},\vec{\mathbf{q}}_n \right)}{\partial\vec{\mathbf{q}}_n} - \left[\frac{\partial \vec{\mathbf{g}}_d (\vec{\mathbf{q}}_n)}{\partial \vec{\mathbf{q}}_n} \right]^\mathrm{T} \vec{\boldsymbol{\uplambda}}_n + \vec{\mathbf{f}}_d^- \left( \vec{\mathbf{q}}_n,\vec{\mathbf{q}}_{n+1} \right) + \vec{\mathbf{f}}_d^+ \left( \vec{\mathbf{q}}_{n-1},\vec{\mathbf{q}}_n \right) &= \vec{\mathbf{0}} \\
\vec{\mathbf{g}}(\vec{\mathbf{q}}_{n+1}) &= \vec{\mathbf{0}}
\end{aligned}
\end{equation}

In order to initiate solving Equations  \ref{eq:eulerLagrange} and \ref{eq:lagrangeDalembert}, both $\vec{\mathbf{q}}_0$ and $\vec{\mathbf{q}}_1$ are required.  In most applications, however, only initial conditions ($\vec{\mathbf{q}}_0$ and $\dot{\vec{\mathbf{q}}}_0$) are provided. To form an equation for $\vec{\mathbf{q}}_1$, the discrete version of the \emph{Legendre transform} can be used. The Legendre transform is a way to switch between the Lagrangian formulation of a dynamic system, which is in the $(\vec{\mathbf{q}}, \dot{\vec{\mathbf{q}}})$ space, to its Hamiltonian formulation in the $(\vec{\mathbf{q}}, \vec{\mathbf{p}})$ space, where $\vec{\mathbf{p}}$ is the system's momentum. If $\vec{\mathbf{p}}_{n,n+1}^{-}$ and $\vec{\mathbf{p}}_{n-1,n}^{+}$ are the so-called \emph{pre}- and \emph{post-momenta} at time $t_n$, for a \emph{constrained dissipative} system, the \emph{discrete Legendre transform} reads as \cite{Leyendecker2008} 

\begin{equation}\label{eq:legendre}
\begin{aligned}
\vec{\mathbf{p}}_{n,n+1}^{-} &= -\frac{\partial L_d \left( \vec{\mathbf{q}}_n,\vec{\mathbf{q}}_{n+1} \right)}{\partial\vec{\mathbf{q}}_n} + \frac{1}{2} \left[\frac{\partial \vec{\mathbf{g}}_d (\vec{\mathbf{q}}_n)}{\partial \vec{\mathbf{q}}_n} \right]^\mathrm{T} \vec{\boldsymbol{\uplambda}}_n - \vec{\mathbf{f}}_d^- \left( \vec{\mathbf{q}}_n,\vec{\mathbf{q}}_{n+1} \right), \\
\vec{\mathbf{p}}_{n-1,n}^{+} &= \frac{\partial L_d \left( \vec{\mathbf{q}}_{n-1},\vec{\mathbf{q}}_{n} \right)}{\partial\vec{\mathbf{q}}_n} - \frac{1}{2} \left[\frac{\partial \vec{\mathbf{g}}_d (\vec{\mathbf{q}}_n)}{\partial \vec{\mathbf{q}}_n} \right]^\mathrm{T} \vec{\boldsymbol{\uplambda}}_n + \vec{\mathbf{f}}_d^+ \left( \vec{\mathbf{q}}_{n-1},\vec{\mathbf{q}}_{n} \right)
\end{aligned}
\end{equation}

The numerical schemes in Equations \ref{eq:lagrangeDalembert} and \ref{eq:eulerLagrange} are regarded as \emph{momentum matching}, meaning that at every $t_n$

\begin{equation}\label{eq:momentum}
\vec{\mathbf{p}}_{n}=\vec{\mathbf{p}}_{n-1,n}^{+}=\vec{\mathbf{p}}_{n,n+1}^{-}
\end{equation}

Thus, as $\vec{\mathbf{p}}_{0}=\mathbf{M}\dot{\vec{\mathbf{q}}}_0$ at $t_0$ and given $\vec{\mathbf{q}}_0$ and $\dot{\vec{\mathbf{q}}}_0$ as inputs, solving the following system of nonlinear equations, with $\mathbf{M}$ as the mass matrix, provides $\vec{\mathbf{q}}_1$ and $\vec{\boldsymbol{\uplambda}}_0$.

\begin{equation}\label{eq:q1lambda0}
\begin{aligned}
\vec{\mathbf{p}}_{0} = -\frac{\partial L_d \left( \vec{\mathbf{q}}_0,\vec{\mathbf{q}}_{1} \right)}{\partial\vec{\mathbf{q}}_0} + \frac{1}{2} \left[\frac{\partial \vec{\mathbf{g}}_d (\vec{\mathbf{q}}_0)}{\partial \vec{\mathbf{q}}_0} \right]^\mathrm{T} \vec{\boldsymbol{\uplambda}}_0 - \vec{\mathbf{f}}_d^- \left( \vec{\mathbf{q}}_0,\vec{\mathbf{q}}_{1} \right) &= \mathbf{M} \dot{\vec{\mathbf{q}}}_0, \\
\vec{\mathbf{g}}(\vec{\mathbf{q}}_1) &= \vec{\mathbf{0}}
\end{aligned}
\end{equation}

The solutions computed for the equations of motions can now be used in sensitivity analysis, which is detailed in the next section.

\section{Sensitivity analysis using the discrete adjoint variable method}
\label{sec:sensitivity}

As mentioned earlier, suppose the optimization goal is to minimize the objective function in Equation \ref{eq:genericform}. Adopting the idea of variational integrators explained in the previous section and assuming a constant time-step size (for simplicity), Equation \ref{eq:genericform} can be discretized using one-point quadrature rule as
\begin{equation}\label{eq:disObj}
\begin{aligned}
\phi &=  F \left( \vec{\mathbf{q}}_0, \dot{\vec{\mathbf{q}}}_0, \vec{\mathbf{q}}_T, \dot{\vec{\mathbf{q}}}_T, \vec{\mathbf{a}} \right) + \displaystyle\sum_{n=0}^{N-1} \int_{t_n}^{t_{n+1}} H \left( \vec{\mathbf{q}}, \dot{\vec{\mathbf{q}}}, \vec{\boldsymbol{\uplambda}}, \vec{\mathbf{a}} \right) dt \\
&= F \left( \vec{\mathbf{q}}_0, \dot{\vec{\mathbf{q}}}_0, \vec{\mathbf{q}}_N, \frac{\vec{\mathbf{q}}_N - \vec{\mathbf{q}}_{N-1}}{h}, \vec{\mathbf{a}} \right) +
\displaystyle\sum_{n=0}^{N-1} \left( hH \left( (1-\alpha)\vec{\mathbf{q}}_n + \alpha \vec{\mathbf{q}}_{n+1}, \frac{\vec{\mathbf{q}}_{n+1}-\vec{\mathbf{q}}_{n}}{h} \right), \vec{\boldsymbol{\uplambda}}_n, \vec{\mathbf{a}} \right) \\
&= \Phi \left( \vec{\mathbf{q}}_0, \vec{\mathbf{q}}_1, \dotsc, \vec{\mathbf{q}}_N, \dot{\vec{\mathbf{q}}}_0, \vec{\boldsymbol{\uplambda}}_0, \vec{\boldsymbol{\uplambda}}_1, \dotsc, \vec{\boldsymbol{\uplambda}}_{N-1}, \vec{\mathbf{a}} \right)
\end{aligned}
\end{equation}
To find a local minimum of this function using a gradient-based approach, its derivative with respect to each design variable is required. For a component $a_i$ of $\vec{\mathbf{a}}$, this is given by
\begin{equation}\label{eq:dObjDai}
\frac{d\Phi}{da_i} = \frac{\partial \Phi}{\partial a_i} + \left[ \frac{\partial \Phi}{\partial \dot{\vec{\mathbf{q}}}_0} \right]^{\mathrm{T}} \frac{\dot{\vec{\mathbf{q}}}_0}{\partial a_i} + \displaystyle\sum_{n=0}^{N} \left( \left[ \frac{\partial \Phi}{\partial \vec{\mathbf{q}}_n} \right]^{\mathrm{T}} \frac{\vec{\mathbf{q}}_n}{\partial a_i} \right) + 
\displaystyle\sum_{n=0}^{N-1} \left( \left[ \frac{\partial \Phi}{\partial \vec{\boldsymbol{\uplambda}}_n} \right]^{\mathrm{T}} \frac{\vec{\boldsymbol{\uplambda}}_n}{\partial a_i} \right)
\end{equation}

Thus, to compute Equation \ref{eq:dObjDai}, one needs to know the derivative of the discrete state variables and Lagrange multipliers with respect to each design variable, namely $\frac{\vec{\mathbf{q}}_n}{\partial a_i} $ and $\frac{\vec{\boldsymbol{\uplambda}}_n}{\partial a_i}$. This can be done through the \emph{discrete adjoint variable method} (DAVM). In this section, the detailed derivation of the equations involved in the DAVM and how to solve them are presented.

Based on Equations \ref{eq:lagrangeDalembert} and \ref{eq:q1lambda0}, the set of motion equations to be solved are of the form

\begin{equation}\label{eq:motionEq}
\begin{aligned}
\vec{\mathbf{c}}_0 \left( \vec{\mathbf{q}}_0, \vec{\mathbf{q}}_1, \vec{\boldsymbol{\uplambda}}_0, \dot{\vec{\mathbf{q}}}_0, \vec{\mathbf{a}}  \right) &= \vec{\mathbf{0}} \\
\vec{\mathbf{g}}_1 = \vec{\mathbf{g}}\left(\vec{\mathbf{q}}_1, \vec{\mathbf{a}} \right) &= \vec{\mathbf{0}} \\
\vec{\mathbf{c}}_n \left( \vec{\mathbf{q}}_{n-1}, \vec{\mathbf{q}}_n, \vec{\mathbf{q}}_{n+1}, \vec{\boldsymbol{\uplambda}}_n, \vec{\mathbf{a}}  \right) &= \vec{\mathbf{0}} \quad \text{for $n=1, 2, \dotsc, N-1$} \\
\vec{\mathbf{g}}_{n+1} = \vec{\mathbf{g}}\left(\vec{\mathbf{q}}_{n+1}, \vec{\mathbf{a}} \right) &= \vec{\mathbf{0}} \quad \text{for $n=1, 2, \dotsc, N-1$} 
\end{aligned}
\end{equation}

Also, at $t=0$, the following equations stand

\begin{equation}\label{eq:g0Dg0}
\begin{aligned}
\vec{\mathbf{g}}_0=\vec{\mathbf{g}}  \left( \vec{\mathbf{q}}_0, \vec{\mathbf{a}} \right) &= \vec{\mathbf{0}} \\
\frac{\partial\vec{\mathbf{g}}  \left( \vec{\mathbf{q}}_0, \vec{\mathbf{a}} \right) }{\partial \vec{\mathbf{q}}_0} \dot{\vec{\mathbf{q}}}_0 &= \vec{\mathbf{0}}
\end{aligned}
\end{equation}

Differentiating each of these six equations with respect to $a_i$ leads to

\begin{equation}\label{eq:DmotionEq}
\begin{aligned}
\frac{\partial \vec{\mathbf{g}}_0}{\partial a_i} + \frac{\partial \vec{\mathbf{g}}_0}{\partial \vec{\mathbf{q}}_0} \frac{\partial \vec{\mathbf{q}}_0}{\partial a_i} &= \vec{\mathbf{0}} \\
\frac{\partial}{\partial a_i} \left(\frac{\partial \vec{\mathbf{g}}_0}{\partial \vec{\mathbf{q}}_0} \right) \dot{\vec{\mathbf{q}}}_0 + 
\frac{\partial \vec{\mathbf{g}}_0}{\partial \vec{\mathbf{q}}_0} \frac{\partial \dot{\vec{\mathbf{q}}}_0}{\partial a_i} + \frac{\partial}{\partial \vec{\mathbf{q}}_0} \left(\frac{\partial \vec{\mathbf{g}}_0}{\partial \vec{\mathbf{q}}_0} \right) \frac{\partial \vec{\mathbf{q}}_0}{\partial a_i} \dot{\vec{\mathbf{q}}}_0 &= \vec{\mathbf{0}} \\
\frac{\partial \vec{\mathbf{c}}_0}{\partial a_i} + \frac{\partial \vec{\mathbf{c}}_0}{\partial \vec{\mathbf{q}}_0} \frac{\partial \vec{\mathbf{q}}_0}{\partial a_i} + 
\frac{\partial \vec{\mathbf{c}}_0}{\partial \vec{\mathbf{q}}_1} \frac{\partial \vec{\mathbf{q}}_1}{\partial a_i} + 
\frac{\partial \vec{\mathbf{c}}_0}{\partial \vec{\boldsymbol{\uplambda}}_0} \frac{\partial \vec{\boldsymbol{\uplambda}}_0}{\partial a_i} + 
\frac{\partial \vec{\mathbf{c}}_0}{\partial \dot{\vec{\mathbf{q}}}_0} \frac{\partial \dot{\vec{\mathbf{q}}}_0}{\partial a_i}
&= \vec{\mathbf{0}} \\
\frac{\partial \vec{\mathbf{g}}_1}{\partial a_i} + \frac{\partial \vec{\mathbf{g}}_1}{\partial \vec{\mathbf{q}}_1} \frac{\partial \vec{\mathbf{q}}_1}{\partial a_i} &= \vec{\mathbf{0}} \\
\frac{\partial \vec{\mathbf{c}}_n}{\partial a_i} + \frac{\partial \vec{\mathbf{c}}_n}{\partial \vec{\mathbf{q}}_{n-1}} \frac{\partial \vec{\mathbf{q}}_{n-1}}{\partial a_i} + 
\frac{\partial \vec{\mathbf{c}}_n}{\partial \vec{\mathbf{q}}_n} \frac{\partial \vec{\mathbf{q}}_n}{\partial a_i} + 
\frac{\partial \vec{\mathbf{c}}_n}{\partial a_i} + \frac{\partial \vec{\mathbf{c}}_n}{\partial \vec{\mathbf{q}}_{n+1}} \frac{\partial \vec{\mathbf{q}}_{n+1}}{\partial a_i} +
\frac{\partial \vec{\mathbf{c}}_n}{\partial \vec{\boldsymbol{\uplambda}}_n} \frac{\partial \vec{\boldsymbol{\uplambda}}_n}{\partial a_i}
&= \vec{\mathbf{0}} \\
\frac{\partial \vec{\mathbf{g}}_{n+1}}{\partial a_i} + \frac{\partial \vec{\mathbf{g}}_{n+1}}{\partial \vec{\mathbf{q}}_{n+1}} \frac{\partial \vec{\mathbf{q}}_{n+1}}{\partial a_i} &= \vec{\mathbf{0}}
\end{aligned} 
\end{equation}

Since the geometrical design variables are considered in this study, the constraint equations $\vec{\mathbf{g}}_n$ are dependent on the design variables, as will be seen in the forthcoming sections. Hence, their derivative with respect to $a_i$ is not zero. 

Note: Equation \ref{eq:DmotionEq} provides a systems of algebraic equations with $\frac{\partial \dot{\vec{\mathbf{q}}}_0}{\partial a_i}$, $\frac{\partial \vec{\mathbf{q}}_n}{\partial a_i}$ and $\frac{\partial \vec{\boldsymbol{\uplambda}}_n}{\partial a_i}$ as unknowns. Having done the forward simulation, by solving this system and implementing the computed values in Equation \ref{eq:dObjDai}, the sensitivity of the objective function can be calculated. This way of performing the sensitivity analysis is referred to as the \emph{discrete direct differentiation method}. The main issue with this technique is that Equation \ref{eq:DmotionEq} needs to be solved separately for all design variables. Therefore, for complex problems with a large number of design variables, it becomes computationally too expensive. For such cases, the DAVM can be utilized, which eliminates the need for computing $\frac{\partial \vec{\mathbf{q}}_n}{\partial a_i}$ and $\frac{\partial \vec{\boldsymbol{\uplambda}}_n}{\partial a_i}$ in Equation \ref{eq:DmotionEq}.

One can introduce \emph{adjoint vectors} $\vec{\boldsymbol{\upmu}}_n \in \mathbb{R}^m \: (n=0,1,\dotsc,N-1)$ and $\vec{\boldsymbol{\upeta}}_n \in \mathbb{R}^l \: (n=1,\dotsc,N)$ associated with the dynamic and constraint equations at each time step {($m$: number of degrees of freedom of the MBS and $l$: number of holonomic constraints)}. Since the expressions in Equation \ref{eq:DmotionEq} are all equal to zero, multiplying them by the transpose of the adjoint vectors and subtracting them from Equation \ref{eq:dObjDai} does not change the value of $\frac{d\Phi}{da_i}$, regardless of the values of $\vec{\boldsymbol{\upmu}}_n$ and $\vec{\boldsymbol{\upeta}}_n$. Thus,
\begin{equation}\label{eq:dPhiDai1}
\begin{aligned}
\frac{d\Phi}{da_i} &= \frac{\partial \Phi}{\partial a_i} + \left[\frac{\partial \Phi}{\partial \dot{\vec{\mathbf{q}}}_0} \right]^\mathrm{T} \frac{\partial \dot{\vec{\mathbf{q}}}_0}{\partial a_i} + \displaystyle\sum_{n=0}^N \left( \left[\frac{\partial \Phi}{\partial \vec{\mathbf{q}}_n} \right]^\mathrm{T} \frac{\partial \vec{\mathbf{q}}_n}{\partial a_i} \right) + \displaystyle\sum_{n=0}^{N-1} \left( \left[\frac{\partial \Phi}{\partial \vec{\boldsymbol{\uplambda}}_n} \right]^\mathrm{T} \frac{\partial \vec{\boldsymbol{\uplambda}}_n}{\partial a_i} \right) \\
&- \displaystyle\sum_{n=1}^{N-1} \left( \vec{\boldsymbol{\upmu}}_n^{\mathrm{T}} \left[ \frac{\partial \vec{\mathbf{c}}_n}{\partial a_i} + \frac{\partial \vec{\mathbf{c}}_n}{\partial \vec{\mathbf{q}}_{n-1}} \frac{\partial \vec{\mathbf{q}}_{n-1}}{\partial a_i} + 
\frac{\partial \vec{\mathbf{c}}_n}{\partial \vec{\mathbf{q}}_n} \frac{\partial \vec{\mathbf{q}}_n}{\partial a_i} + 
\frac{\partial \vec{\mathbf{c}}_n}{\partial a_i} + \frac{\partial \vec{\mathbf{c}}_n}{\partial \vec{\mathbf{q}}_{n+1}} \frac{\partial \vec{\mathbf{q}}_{n+1}}{\partial a_i} +
\frac{\partial \vec{\mathbf{c}}_n}{\partial \vec{\boldsymbol{\uplambda}}_n} \frac{\partial \vec{\boldsymbol{\uplambda}}_n}{\partial a_i} \right] \right) \\
&- \displaystyle\sum_{n=1}^{N} \left( \vec{\boldsymbol{\upeta}}_n^{\mathrm{T}} \left[\frac{\partial \vec{\mathbf{g}}_{n}}{\partial a_i} + \frac{\partial \vec{\mathbf{g}}_{n}}{\partial \vec{\mathbf{q}}_{n}} \frac{\partial \vec{\mathbf{q}}_{n}}{\partial a_i} \right] \right) \\
&- \vec{\boldsymbol{\upmu}}_0^\mathrm{T} \left[ \frac{\partial \vec{\mathbf{c}}_0}{\partial a_i} + \frac{\partial \vec{\mathbf{c}}_0}{\partial \vec{\mathbf{q}}_0} \frac{\partial \vec{\mathbf{q}}_0}{\partial a_i} + 
\frac{\partial \vec{\mathbf{c}}_0}{\partial \vec{\mathbf{q}}_1} \frac{\partial \vec{\mathbf{q}}_1}{\partial a_i} + 
\frac{\partial \vec{\mathbf{c}}_0}{\partial \vec{\boldsymbol{\uplambda}}_0} \frac{\partial \vec{\boldsymbol{\uplambda}}_0}{\partial a_i} + 
\frac{\partial \vec{\mathbf{c}}_0}{\partial \dot{\vec{\mathbf{q}}}_0} \frac{\partial \dot{\vec{\mathbf{q}}}_0}{\partial a_i} \right]
\end{aligned} 
\end{equation}
Re-indexing the terms and putting the similar ones together lead to
\begin{equation}\label{eq:dPhiDai2}
\begin{aligned}
\frac{d\Phi}{da_i} &= \frac{\partial \Phi}{\partial a_i} - \displaystyle\sum_{n=0}^{N-1} \left( \vec{\boldsymbol{\upmu}}_n^\mathrm{T} \frac{\partial \vec{\mathbf{c}}_n}{\partial a_i} \right) - \displaystyle\sum_{n=1}^{N} \left( \vec{\boldsymbol{\upeta}}_n^\mathrm{T} \frac{\partial \vec{\mathbf{g}}_n}{\partial a_i} \right) \\
&+ \displaystyle\sum_{n=1}^{N-2} \left( \left[ \left[\frac{\partial \Phi}{\partial \vec{\mathbf{q}}_n} \right]^\mathrm{T} - \vec{\boldsymbol{\upmu}}_{n+1}^{\mathrm{T}} \frac{\partial \vec{\mathbf{c}}_{n+1}}{\partial \vec{\mathbf{q}}_n} - \vec{\boldsymbol{\upmu}}_{n}^{\mathrm{T}} \frac{\partial \vec{\mathbf{c}}_{n}}{\partial \vec{\mathbf{q}}_n} - \vec{\boldsymbol{\upmu}}_{n-1}^{\mathrm{T}} \frac{\partial \vec{\mathbf{c}}_{n-1}}{\partial \vec{\mathbf{q}}_n} - \vec{\boldsymbol{\upeta}}_{n}^{\mathrm{T}} \frac{\partial \vec{\mathbf{g}}_{n}}{\partial \vec{\mathbf{q}}_n} \right] \frac{\partial \vec{\mathbf{q}}_n}{\partial a_i} \right) \\
&+ \displaystyle\sum_{n=0}^{N-1} \left( \left[ \left[\frac{\partial \Phi}{\partial \vec{\boldsymbol{\uplambda}}_n} \right]^\mathrm{T} - \vec{\boldsymbol{\upmu}}_{n}^{\mathrm{T}} \frac{\partial \vec{\mathbf{c}}_{n}}{\partial \vec{\boldsymbol{\uplambda}}_n} \right] \frac{\partial \vec{\boldsymbol{\uplambda}}_n}{\partial a_i} \right) \\
&+ \left[ \left[\frac{\partial \Phi}{\partial \vec{\mathbf{q}}_{N-1}} \right]^\mathrm{T} - \vec{\boldsymbol{\upmu}}_{N-1}^{\mathrm{T}} \frac{\partial \vec{\mathbf{c}}_{N-1}}{\partial \vec{\mathbf{q}}_{N-1}} - \vec{\boldsymbol{\upmu}}_{N-2}^{\mathrm{T}} \frac{\partial \vec{\mathbf{c}}_{N-2}}{\partial \vec{\mathbf{q}}_{N-1}} - \vec{\boldsymbol{\upeta}}_{N-1}^{\mathrm{T}} \frac{\partial \vec{\mathbf{g}}_{N-1}}{\partial \vec{\mathbf{q}}_{N-1}} \right] \frac{\partial \vec{\mathbf{q}}_{N-1}}{\partial a_i} \\
&+ \left[ \left[\frac{\partial \Phi}{\partial \vec{\mathbf{q}}_{N}} \right]^\mathrm{T} - \vec{\boldsymbol{\upmu}}_{N-1}^{\mathrm{T}} \frac{\partial \vec{\mathbf{c}}_{N-1}}{\partial \vec{\mathbf{q}}_{N}} - \vec{\boldsymbol{\upeta}}_{N}^{\mathrm{T}} \frac{\partial \vec{\mathbf{g}}_{N}}{\partial \vec{\mathbf{q}}_{N}} \right] \frac{\partial \vec{\mathbf{q}}_{N}}{\partial a_i} \\
&+ \left[ \left[\frac{\partial \Phi}{\partial \vec{\mathbf{q}}_{0}} \right]^\mathrm{T} - \vec{\boldsymbol{\upmu}}_{1}^{\mathrm{T}} \frac{\partial \vec{\mathbf{c}}_{1}}{\partial \vec{\mathbf{q}}_{0}} - \vec{\boldsymbol{\upmu}}_{0}^{\mathrm{T}} \frac{\partial \vec{\mathbf{c}}_{0}}{\partial \vec{\mathbf{q}}_{0}} \right] \frac{\partial \vec{\mathbf{q}}_{0}}{\partial a_i} \\ 
&+ \left[ \left[\frac{\partial \Phi}{\partial \dot{\vec{\mathbf{q}}}_0} \right]^\mathrm{T} - \vec{\boldsymbol{\upmu}}_{0}^{\mathrm{T}} \frac{\partial \vec{\mathbf{c}}_{0}}{\partial \dot{\vec{\mathbf{q}}}_0} \right] \frac{\partial \dot{\vec{\mathbf{q}}}_0}{\partial a_i}
\end{aligned} 
\end{equation}
As this equation holds for any values of $\vec{\boldsymbol{\upmu}}_n$ and $\vec{\boldsymbol{\upeta}}_n$, it is possible to choose these in such a way that bypasses the computation of the gradients of the state variables and Lagrange multipliers in the sensitivity analysis.  That is, one can choose $\vec{\boldsymbol{\upmu}}_n$ and $\vec{\boldsymbol{\upeta}}_n$ so that the coefficients of $\frac{\partial \vec{\mathbf{q}}_n}{\partial a_i}$ and $\frac{\partial \vec{\boldsymbol{\uplambda}}_n}{\partial a_i}$ in Equation \ref{eq:dPhiDai2} all become zero. To do so, first, $\vec{\boldsymbol{\upmu}}_{N-1}$, $\vec{\boldsymbol{\upmu}}_{N-2}$, $\vec{\boldsymbol{\upeta}}_{N}$ and $\vec{\boldsymbol{\upeta}}_{N-1}$ are found from the following two systems of linear algebraic equations

\begin{equation}\label{eq:adjointsystem1}
\begin{aligned}
& \begin{dcases}
\left[ \frac{\partial \vec{\mathbf{c}}_{N-1}}{\partial \vec{\mathbf{q}}_{N}} \right]^{\mathrm{T}} \vec{\boldsymbol{\upmu}}_{N-1}+ \left[ \frac{\partial \vec{\mathbf{g}}_{N}}{\partial \vec{\mathbf{q}}_{N}} \right]^{\mathrm{T}} \vec{\boldsymbol{\upeta}}_{N} \hspace{2.8em} = \frac{\partial \Phi}{\partial \vec{\mathbf{q}}_{N}} \\[0.5em]
\left[ \frac{\partial \vec{\mathbf{c}}_{N-1}}{\partial \vec{\boldsymbol{\uplambda}}_{N-1}} \right]^{\mathrm{T}} \vec{\boldsymbol{\upmu}}_{N-1} \hspace{9.2em} = \frac{\partial \Phi}{\partial \vec{\boldsymbol{\uplambda}}_{N-1}} \\[0.5em]
\end{dcases} \\[0.5em]
& \begin{dcases}
\left[ \frac{\partial \vec{\mathbf{c}}_{N-2}}{\partial \vec{\mathbf{q}}_{N-1}} \right]^{\mathrm{T}} \vec{\boldsymbol{\upmu}}_{N-2}+ \left[ \frac{\partial \vec{\mathbf{g}}_{N-1}}{\partial \vec{\mathbf{q}}_{N-1}} \right]^{\mathrm{T}} \vec{\boldsymbol{\upeta}}_{N-1} &= \frac{\partial \Phi}{\partial \vec{\mathbf{q}}_{N-1}} - \left[ \frac{\partial \vec{\mathbf{c}}_{N-1}}{\partial \vec{\mathbf{q}}_{N-1}} \right]^{\mathrm{T}} \vec{\boldsymbol{\upmu}}_{N-1} \\[0.5em]
\left[ \frac{\partial \vec{\mathbf{c}}_{N-2}}{\partial \vec{\boldsymbol{\uplambda}}_{N-2}} \right]^{\mathrm{T}} \vec{\boldsymbol{\upmu}}_{N-2} &= \frac{\partial \Phi}{\partial \vec{\boldsymbol{\uplambda}}_{N-2}} \\[0.5em]
\end{dcases}
\end{aligned}
\end{equation}

Then, $\vec{\boldsymbol{\upmu}}_{n-1}$ and $\vec{\boldsymbol{\upeta}}_{n} (n=N-2, N-3, \dotsc, 1)$ can be computed by solving the system of linear algebraic equations below 

\begin{equation}\label{eq:adjointsystem2}
\begin{dcases}
\left[ \frac{\partial \vec{\mathbf{c}}_{n-1}}{\partial \vec{\mathbf{q}}_n} \right]^{\mathrm{T}} \vec{\boldsymbol{\upmu}}_{n-1} + \left[ \frac{\partial \vec{\mathbf{g}}_{n}}{\partial \vec{\mathbf{q}}_n} \right]^{\mathrm{T}} \vec{\boldsymbol{\upeta}}_{n} &= \frac{\partial \Phi}{\partial \vec{\mathbf{q}}_n} - \left[ \frac{\partial \vec{\mathbf{c}}_{n+1}}{\partial \vec{\mathbf{q}}_n} \right]^{\mathrm{T}} \vec{\boldsymbol{\upmu}}_{n+1} - \left[ \frac{\partial \vec{\mathbf{c}}_{n}}{\partial \vec{\mathbf{q}}_n} \right]^{\mathrm{T}} \vec{\boldsymbol{\upmu}}_{n} \\
\left[ \frac{\partial \vec{\mathbf{c}}_{n-1}}{\partial \vec{\boldsymbol{\uplambda}}_{n-1}} \right]^{\mathrm{T}} \vec{\boldsymbol{\upmu}}_{n-1} &= \frac{\partial \Phi}{\partial \vec{\boldsymbol{\uplambda}}_{n-1}} \\
\end{dcases}
\end{equation}

Note: Equations \ref{eq:adjointsystem1} and \ref{eq:adjointsystem2} are solved in a backward manner from $n=N$ to $n=1$ with $\vec{\boldsymbol{\upmu}}_{n-1}$ and $\vec{\boldsymbol{\upeta}}_n$ as unknowns. That is why the process of solving the adjoint equations is sometimes called \emph{backward simulation}. Unlike Equation \ref{eq:DmotionEq} in the discrete direct differentiation method, systems of equations in Equations \ref{eq:adjointsystem1} and \ref{eq:adjointsystem2} need to be solved only once for each iteration of optimization. Once $\vec{\boldsymbol{\upmu}}_n$ and $\vec{\boldsymbol{\upeta}}_n$ are computed and $\frac{\partial \vec{\mathbf{q}}_0}{\partial a_i}$ and $\frac{\dot{\vec{\mathbf{q}}}_0}{\partial a_i}$ are calculated using the first two expressions in Equation \ref{eq:DmotionEq}, the sensitivity Equation \ref{eq:dPhiDai2} reduces to
\begin{equation}\label{eq:finalSensitivity}
\begin{aligned}
\frac{d\Phi}{da_i} &= \frac{\partial \Phi}{\partial a_i} - \displaystyle\sum_{n=0}^{N-1} \left( \vec{\boldsymbol{\upmu}}_n^\mathrm{T} \frac{\partial \vec{\mathbf{c}}_n}{\partial a_i} \right) - \displaystyle\sum_{n=1}^{N} \left( \vec{\boldsymbol{\upeta}}_n^\mathrm{T} \frac{\partial \vec{\mathbf{g}}_n}{\partial a_i} \right) \\
&+ \left[ \left[\frac{\partial \Phi}{\partial \vec{\mathbf{q}}_{0}} \right]^\mathrm{T} - \vec{\boldsymbol{\upmu}}_{1}^{\mathrm{T}} \frac{\partial \vec{\mathbf{c}}_{1}}{\partial \vec{\mathbf{q}}_{0}} - \vec{\boldsymbol{\upmu}}_{0}^{\mathrm{T}} \frac{\partial \vec{\mathbf{c}}_{0}}{\partial \vec{\mathbf{q}}_{0}} \right] \frac{\partial \vec{\mathbf{q}}_{0}}{\partial a_i} \\ 
&+ \left[ \left[\frac{\partial \Phi}{\partial \dot{\vec{\mathbf{q}}}_0} \right]^\mathrm{T} - \vec{\boldsymbol{\upmu}}_{0}^{\mathrm{T}} \frac{\partial \vec{\mathbf{c}}_{0}}{\partial \dot{\vec{\mathbf{q}}}_0} \right] \frac{\partial \dot{\vec{\mathbf{q}}}_0}{\partial a_i}
\end{aligned} 
\end{equation}
This is the final equation for the derivative of the objective function with respect to design variable $a_i$. The same process can be applied to compute the gradients of optimization constraint functions.

\section{Sensitivity with respect to geometrical design variables}
\label{sec:SensitivityGeo}

As aforementioned, geometrical design variables are those affecting the shape of MBS configurations, such as the dimensions of components and joint positions. Unlike non-geometrical parameters whose values influence only one component of an assembly, geometrical parameters could be related to multiple bodies that share that variable. For instance, changing the position of a joint would impact all the components connected to it. For the sake of brevity and to narrow the premise, this section focuses on MBSs composed of springs, dampers, beams and rigid bodies interacting via spherical (pin) and welded (fixed) joints. The same procedure can be extended to handle more complicated bodies and joints. In this section, first, a brief overview of \emph{rotation-free} formulations for beam and rigid body dynamics, constraint equations between them and the equations of springs and dampers is provided. Then, as an example, the sensitivity analysis of a simple rigid-flexible multibody considering geometrical design variables is described. The avoidance of rotation parameters leads to a \emph{constant} mass matrix for each body. This turns out to be highly beneficial toward having a conserving numerical time-stepping solver \cite{Shabana2001, Uhlar2009}, like the one used in this paper, but also toward simplifying the equations of motion, thus facilitating the required gradient computations with respect to geometrical and non-geometrical design variables. 

\subsection{Absolute nodal coordinate formulation for beams}
\label{subsec:absolute}

Initially proposed by Shabana \cite{Shabana1996}, the  \emph{absolute nodal coordinate formulation} (ANCF) uses the \emph{positions} and \emph{slopes} of nodes in the global inertial frame as the generalized coordinates for beams. This results in a constant mass matrix for beams, and subsequently cancels the nonlinear terms of centrifugal and Coriolis inertia forces in the equations of motion. It also leads to simple expressions for the constraint equations, which is of a great benefit for defining different types of joints and including geometrical design variables in the optimization routine. The ANCF has been widely implemented in modeling deformable objects, such as beams, plates and solids, in dynamic problems. It is capable of describing the rigid body modes and accurately solving large deformation problems \cite{Gerstmayr2006, Garcia-Vallejo2003, Gerstmayr2013, Shabana2001, Yakoub2001}. In this paper, the basics of this theory for a two-noded gradient-deficient Euler-Bernoulli beam element with uniform cross-section is provided \cite{Berzeri2000}. 

\begin{figure}
\centering
\includegraphics[width=0.45\textwidth]{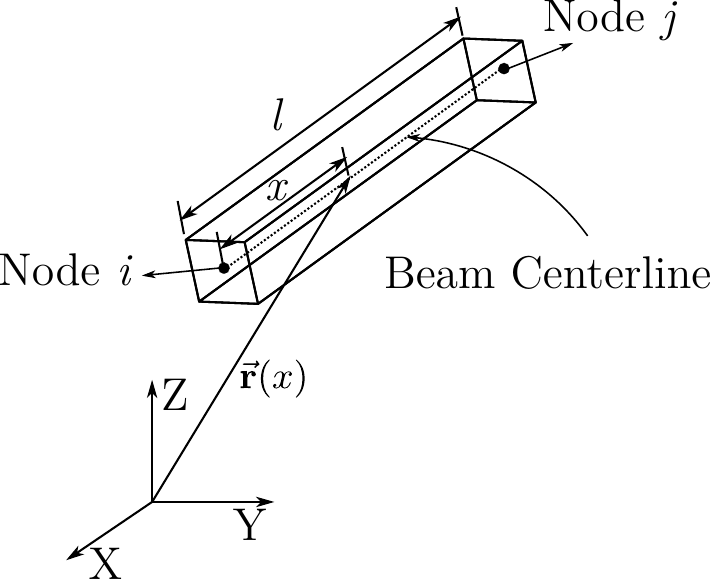}
\caption{\label{fig:beam}Two-noded Euler-Bernoulli beam element}
\end{figure}

Consider a beam element as depicted in Figure \ref{fig:beam}. For this element, the position vector $\vec{\mathbf{r}}$ of an arbitrary point in the \emph{global} Cartesian frame can be written as
\begin{equation}\label{eq:position}
\vec{\mathbf{r}}(x,t)=\mathbf{S}(x)\vec{\mathbf{q}}(t)
\end{equation}

where $x$ is the local coordinate of that point along the center-line in the \emph{undeformed} configuration. In Equation \ref{eq:position}, $\vec{\mathbf{q}}(t)$ is the \emph{time-dependent} vector of nodal degrees of freedom and $\mathbf{S}$ is the \emph{time-independent} shape function matrix expressed via 
\begin{equation}\label{eq:qs}
\begin{aligned}
\vec{\mathbf{q}}(t) &= \begin{bmatrix}
\vec{\mathbf{r}}_i^\mathrm{T} & \frac{\partial\vec{\mathbf{r}}_i^\mathrm{T}}{\partial x} & \vec{\mathbf{r}}_j^\mathrm{T} & \frac{\partial\vec{\mathbf{r}}_j^\mathrm{T}}{\partial x}
\end{bmatrix}^\mathrm{T} \\
\mathbf{S}(x) &= \begin{bmatrix}
s_1 \mathbf{I}_{3\times 3} & s_2 \mathbf{I}_{3\times 3} & s_3 \mathbf{I}_{3\times 3} & s_4 \mathbf{I}_{3\times 3}
\end{bmatrix}\\
s_1 &= 1-3\xi^2+2\xi^3 \\
s_2 &= l(\xi-2\xi^2+\xi^3) \\
s_3 &= 3\xi^2-2\xi^3 \\
s_4 &= l(\xi^3-\xi^2)
\end{aligned}
\end{equation}
in which $\mathbf{I}_{3\times3} $ is the identity matrix, $l$ is the element’s length in the undeformed configuration and $\xi=x/l $. To solve the equations of motion described in Section \ref{sec:Dynamicsof}, the element’s mass matrix and potential energy are required. The former reads as {(details are provided in Appendix {\ref{app:beammass}})}

\begin{equation}\label{eq:kineticdeltaU}
\begin{gathered}
\mathbf{M}_{\text{element}} = \frac{\rho A l}{420} \begin{bmatrix}
156 \mathbf{I}_{3 \times 3} & 22l \mathbf{I}_{3 \times 3} & 54 \mathbf{I}_{3 \times 3} & -13l \mathbf{I}_{3 \times 3} \\[0.8em]
22l \mathbf{I}_{3 \times 3} & 4l^2 \mathbf{I}_{3 \times 3} & 13l \mathbf{I}_{3 \times 3} & -3l^2 \mathbf{I}_{3 \times 3} \\[0.8em]
54 \mathbf{I}_{3 \times 3} & 13l \mathbf{I}_{3 \times 3} & 156 \mathbf{I}_{3 \times 3} & -22l \mathbf{I}_{3 \times 3} \\[0.8em]
-13l \mathbf{I}_{3 \times 3} & -3l^2 \mathbf{I}_{3 \times 3} & -22l \mathbf{I}_{3 \times 3} & 4l^2 \mathbf{I}_{3 \times 3} \\
\end{bmatrix}
\end{gathered}
\end{equation}

where $\rho$, $V$ and $A$ are the density, volume and cross-sectional area of the element, respectively. The potential energy of the element is due to its elastic deformation and gravity. It is given by

\begin{equation}\label{eq:deltaU}
U_{\text{element}} = U_{\text{longitudinal}} + U_{\text{transverse}} + U_{\text{gravity}}= \frac{1}{2}\int \limits_V \left( EA\varepsilon^2+E I\kappa^2 \right) dV + \int\limits_V -\rho \vec{\boldsymbol{\upnu}}^{\mathrm{T}}\vec{\mathbf{r}} \: dV 
\end{equation}

in which $E$, $I$ and $\vec{\boldsymbol{\upnu}}$ are, respectively, Young's modulus, the second moment of area and the gravity vector. In Equation \ref{eq:deltaU}, $\varepsilon$ and $\kappa$ are the element’s longitudinal strain and spatial curvature formulated as

\begin{equation}\label{eq:ek1}
\begin{aligned}
\varepsilon &= \frac{1}{2} \left( \frac{\partial \vec{\mathbf{r}}^\mathrm{T}}{\partial x} \frac{\partial \vec{\mathbf{r}}}{\partial x} -1 \right), \quad \kappa = \frac{\left\| \frac{\partial \vec{\mathbf{r}}}{\partial x} \times \frac{\partial^2 \vec{\mathbf{r}}}{\partial x^2} \right\|}{\left\| \frac{\partial \vec{\mathbf{r}}}{\partial x} \right\|^3} 
\end{aligned}
\end{equation}

In case of small axial deformation and assuming constant longitudinal strain, both expressions in Equation \ref{eq:ek1} can be simplified extensively as follows

\begin{equation}\label{eq:ek2}
\begin{aligned}
\varepsilon &= \frac{d}{l}-1 = \frac{1}{l} \left\|\vec{\mathbf{r}}_j - \vec{\mathbf{r}}_i \right\| -1, \quad \kappa = \left\| \frac{\partial^2 \vec{\mathbf{r}}}{\partial x^2} \right\|
\end{aligned}
\end{equation}

To perform the forward simulation, as in Section \ref{sec:Dynamicsof}, one needs to compute the derivative of element's potential energy, Equation \ref{eq:deltaU}, with respect to $\vec{\mathbf{q}}$, also called the vector of \emph{generalized forces}. Different assumptions and formulas have been proposed for this quantity, some of which could be found in \cite{Berzeri2000, Takahashi1999}. Assuming small axial deformations and using Equations \ref{eq:position}, \ref{eq:deltaU} and \ref{eq:ek2}, the generalized elastic forces take a linear form as 

\begin{equation}\label{eq:linearForm}
\begin{aligned}
\frac{\partial U_{\text{elastic}}}{\partial \vec{\mathbf{q}}} &= \frac{\partial U_{\text{longitudinal}}}{\partial \vec{\mathbf{q}}} + \frac{\partial U_{\text{transverse}}}{\partial \vec{\mathbf{q}}} = \mathbf{K}_{\text{longitudinal}}\vec{\mathbf{q}} + \mathbf{K}_{\text{transverse}}\vec{\mathbf{q}} \\
\end{aligned}
\end{equation}

Expressions of $\mathbf{K}_{\text{longitudinal}}$ and $\mathbf{K}_{\text{transverse}}$ are provided in {Appendix {\ref{app:beammass}}}. The generalized gravity force vector can be computed by

\begin{equation}\label{eq:gengravbeam}
\begin{aligned}
\frac{\partial U_{\text{gravity}}}{\partial \vec{\mathbf{q}}} &= - \left[ \int\limits_V \left[\rho \vec{\boldsymbol{\upnu}}^{\mathrm{T}}\mathbf{S} \right]dV \right] = - \frac{1}{12} \rho A l \begin{bmatrix}
6\mathbf{I}_{3 \times 3} & l\mathbf{I}_{3 \times 3} & 6 \mathbf{I}_{3 \times 3} & -l\mathbf{I}_{3 \times 3}
\end{bmatrix}^{\mathrm{T}} \vec{\boldsymbol{\upnu}}
\end{aligned}
\end{equation}

\subsection{Natural coordinates for rigid bodies}
\label{subsec:natural}

\emph{Natural coordinates}, also known as \emph{basic coordinates} or \emph{Cartesian coordinates}, achieve the rotation-free expression for dynamic rigid bodies by introducing a set of redundant degress of freedom and additional constrains into the equations of motion \cite{DeJalon1987, DeJalon2012, Betsch2001, DeJalon2007}. Similar to ANCF for beams, natural coordinates formulation also results in a constant mass matrix for rigid bodies and extensively simplifies the dynamic and sensitivity equations. 

\begin{figure}
\centering
\includegraphics[width=0.33\textwidth]{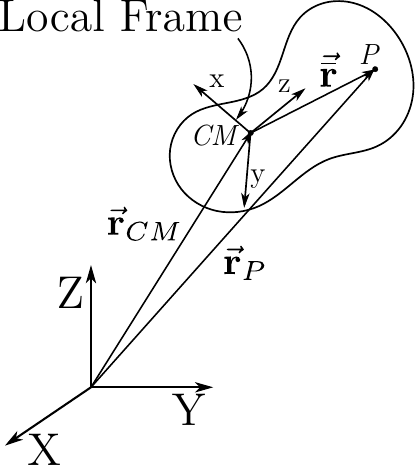}
\caption{\label{fig:rigid}Rigid body}
\end{figure}

Let $\vec{\mathbf{r}}_{CM}$ and $\mathbf{R}$ denote the global position of the {{body’s} center of mass (CM)} and the local frame attached to it, respectively, as shown in Figure \ref{fig:rigid}. If $\vec{\bar{\mathbf{r}}}= \begin{bmatrix}\bar{x} & \bar{y} & \bar{z}\end{bmatrix}^\mathrm{T}$ is the local position vector of an arbitrary point P on the body, its global position vector is given by

\begin{equation}\label{eq:rp}
\vec{\mathbf{r}}_P = \vec{\mathbf{r}}_{CM} + \mathbf{R}\vec{\bar{\mathbf{r}}} = \vec{\mathbf{r}}_{CM} + \mathbf{R} \begin{Bmatrix} \bar{x} \\ \bar{y} \\ \bar{z} \end{Bmatrix}
\end{equation}

Suppose $\vec{\mathbf{e}}_1$, $\vec{\mathbf{e}}_2$ and $\vec{\mathbf{e}}_3$ are the unit basis vectors of the local frame. Then, Equation \ref{eq:rp} can be rewritten as

\begin{equation}\label{eq:rp2}
\vec{\mathbf{r}}_P \left(\bar{x},\bar{y},\bar{z} , t\right)= \begin{bmatrix}
\mathbf{I}_{3\times 3} & \bar{x} \mathbf{I}_{3\times 3} & \bar{y} \mathbf{I}_{3\times 3} & \bar{z} \mathbf{I}_{3\times 3}
\end{bmatrix} \begin{Bmatrix}
\vec{\mathbf{r}}_{CM} \\ \vec{\mathbf{e}}_1 \\ \vec{\mathbf{e}}_2 \\ \vec{\mathbf{e}}_3
\end{Bmatrix} = \mathbf{S}\left( \bar{x},\bar{y},\bar{z} \right)\vec{\mathbf{q}}(t)
\end{equation}

where $\vec{\mathbf{q}}$ is the \emph{time-dependent} vector of the rigid body's degrees of freedom (generalized coordinates). As only 6 degrees of freedom are enough to fully configure a rigid body in the 3D space and $\vec{\mathbf{q}}$ has 12 elements, a set of additional constraints, called \emph{internal constraints}, are required. They can be stated as 

\begin{equation}\label{eq:ginternal}
\vec{\mathbf{g}}_{\text{internal}}(\vec{\mathbf{q}}) = \begin{Bmatrix}
\vec{\mathbf{e}}_1\cdot\vec{\mathbf{e}}_1-1 \\ \vec{\mathbf{e}}_2\cdot\vec{\mathbf{e}}_2-1 \\ \vec{\mathbf{e}}_3\cdot\vec{\mathbf{e}}_3-1 \\ \vec{\mathbf{e}}_1\cdot\vec{\mathbf{e}}_2 \\ \vec{\mathbf{e}}_1\cdot\vec{\mathbf{e}}_3 \\ \vec{\mathbf{e}}_2\cdot\vec{\mathbf{e}}_3
\end{Bmatrix} = \vec{\mathbf{0}}_6
\end{equation}

Internal constraints, in fact, enforce the \emph{orthonormality} of the rigid body's local coordinate system. If the local frame's axes are aligned with the body's principal axes, {incorporating the procedure described in Appendix {\ref{app:rigid}}}, the rigid body's mass matrix takes the following diagonal form

\begin{equation}\label{eq:massRigid}
\mathbf{M}_{\text{rigid}} = \begin{bmatrix}
m\mathbf{I}_{3 \times 3} & \mathbf{0}_{3 \times 3} & \mathbf{0}_{3 \times 3} & \mathbf{0}_{3 \times 3} \\[0.8em]
\mathbf{0}_{3 \times 3} & \frac{1}{2}(I_2+I_3-I_1) \mathbf{I}_{3 \times 3} & \mathbf{0}_{3 \times 3} & \mathbf{0}_{3 \times 3} \\[0.8em]
\mathbf{0}_{3 \times 3} & \mathbf{0}_{3 \times 3} & \frac{1}{2}(I_1+I_3-I_2) \mathbf{I}_{3 \times 3} & \mathbf{0}_{3 \times 3} \\[0.8em]
\mathbf{0}_{3 \times 3} & \mathbf{0}_{3 \times 3} & \mathbf{0}_{3 \times 3} & \frac{1}{2}(I_1+I_2-I_3) \mathbf{I}_{3 \times 3} \\
\end{bmatrix}
\end{equation}

in which, $m$ and $I_i  (i=1,2,3)$ are the body’s mass and principal moments of inertia, respectively. Since there is no time-term in the mass matrix equation, it remains constant during the forward simulation. To solve the equations described in Section \ref{sec:Dynamicsof}, the derivative of body's potential energy, which is only due to gravity, with respect to $\vec{\mathbf{q}}$ is needed. This parameter is called the \emph{generalized gravity force vector} and reads as

\begin{equation}
\begin{aligned}
\frac{\partial U_{\text{gravity}}}{\partial \vec{\mathbf{q}}} &= - \left[ \int\limits_V \left[\rho \vec{\boldsymbol{\upnu}}^{\mathrm{T}}\mathbf{S} \right]dV \right]= - m \begin{bmatrix}
\upnu_x & \upnu_y & \upnu_z & 0 & 0 & 0 & 0 & 0 & 0 & 0 & 0 & 0
\end{bmatrix}^{\mathrm{T}}
\end{aligned}
\end{equation}

where $\vec{\boldsymbol{\upnu}}=\begin{bmatrix}\upnu_x & \upnu_y & \upnu_z \end{bmatrix}^{\mathrm{T}}$ is the gravity vector,

\subsection{Constraint equations using ANCF and natural coordinates}
\label{subsec:Constraint}

A key benefit of utilizing these two rotation-free theories is the simplification they make on the constraint equations. For example, a spherical joint between a point $\mathrm{A}$ on a rigid body and an end point $\mathrm{B}$ on a beam is expressed as

\begin{equation}\label{eq:gSpherical}
\begin{aligned}
\vec{\mathbf{g}}_{\text{spherical}} &= \vec{\mathbf{r}}_{\mathrm{A}} - \vec{\mathbf{r}}_{\mathrm{B}} = \left. \mathbf{S}_{\text{rigid}} \right|_{\mathrm{A}} \vec{\mathbf{q}}_{\text{rigid}} - \left. \mathbf{S}_{\text{beam}} \right|_{\mathrm{B}} \vec{\mathbf{q}}_{\text{beam}}  \\ 
&= \begin{bmatrix}
\left. \mathbf{S}_{\text{rigid}} \right|_{\mathrm{A}} & -\left. \mathbf{S}_{\text{beam}} \right|_{\mathrm{B}}
\end{bmatrix} \begin{Bmatrix}
\vec{\mathbf{q}}_{\text{rigid}} \\ \vec{\mathbf{q}}_{\text{beam}}
\end{Bmatrix} = \vec{\mathbf{0}}_3
\end{aligned}
\end{equation}

and if instead, this connection is welded, Equation \ref{eq:gSpherical} is augmented by these three more equations

\begin{equation}\label{eq:gWelded}
\vec{\mathbf{g}}_{\text{welded}} = \left. \mathbf{R}_{\text{rigid}}^\mathrm{T} \right|_{t_{\text{current}}} \left. \frac{\partial \vec{\mathbf{r}}_{\mathrm{B}}}{\partial x} \right|_{t_{\text{current}}} -
\left. \mathbf{R}_{\text{rigid}}^\mathrm{T} \right|_{t_{0}} \left. \frac{\partial \vec{\mathbf{r}}_{\mathrm{B}}}{\partial x} \right|_{t_{0}} = \vec{\mathbf{0}}_3
\end{equation}

Other types of constraints can be defined in the same fashion \cite{Uhlar2009, Garcia-Vallejo2008}. Using ANCF and natural coordinates, the constraint Jacobian matrix can be analytically evaluated in a straightforward manner.

\subsection{Spring and damper equations}
\label{subsec:springdamper}

Assume a spring and a damper, between a point C of a rigid body and an end point D of a beam. If $k$ is the spring constant and $c$ denotes the damping coefficient, the \emph{generalized spring} and \emph{damper force} vectors are expressed via

\begin{equation}\label{eq:genspring}
\begin{aligned}
\frac{\partial U_{\text{spring}}}{\partial \vec{\mathbf{q}}} =k \frac{l_{t_{\mathrm{current}}}-l_0}{l_{t_{\mathrm{current}}}} &\begin{bmatrix} \mathbf{S}_{\text{rigid}}^{\mathrm{T}} \Big|_{\mathrm{C}}  \mathbf{S}_{\text{rigid}} \Big|_{\mathrm{C}} & -\mathbf{S}_{\text{rigid}}^{\mathrm{T}} \Big|_{\mathrm{C}} \mathbf{S}_{\text{beam}} \Big|_{\mathrm{D}} \\[1.0em]
-\mathbf{S}_{\text{beam}}^{\mathrm{T}} \Big|_{\mathrm{D}}  \mathbf{S}_{\text{rigid}} \Big|_{\mathrm{C}} & \mathbf{S}_{\text{beam}}^{\mathrm{T}} \Big|_{\mathrm{D}} \mathbf{S}_{\text{beam}} \Big|_{\mathrm{D}}
\end{bmatrix} \begin{Bmatrix}
\vec{\mathbf{q}}_{\text{rigid}} \\[1.0em] \vec{\mathbf{q}}_{\text{beam}}
\end{Bmatrix}\\
\vec{\mathbf{f}}_{\text{damper}} = -c &\begin{bmatrix} \mathbf{S}_{\text{rigid}}^{\mathrm{T}} \Big|_{\mathrm{C}}  \mathbf{S}_{\text{rigid}} \Big|_{\mathrm{C}} & -\mathbf{S}_{\text{rigid}}^{\mathrm{T}} \Big|_{\mathrm{C}} \mathbf{S}_{\text{beam}} \Big|_{\mathrm{D}} \\[1.0em]
-\mathbf{S}_{\text{beam}}^{\mathrm{T}} \Big|_{\mathrm{D}}  \mathbf{S}_{\text{rigid}} \Big|_{\mathrm{C}} & \mathbf{S}_{\text{beam}}^{\mathrm{T}} \Big|_{\mathrm{D}} \mathbf{S}_{\text{beam}} \Big|_{\mathrm{D}}
\end{bmatrix} \begin{Bmatrix}
\dot{\vec{\mathbf{q}}}_{\text{rigid}} \\[1.0em] \dot{\vec{\mathbf{q}}}_{\text{beam}}
\end{Bmatrix}
\end{aligned}
\end{equation}

where $l_{t_{\mathrm{current}}}$ and $l_0$ are the current and initial lengths of the spring, respectively. {The derivation of Equation {\ref{eq:genspring}} is provided in Appendix {\ref{app:spring}}.}

\subsection{Example: sensitivity analysis for a simple rigid-spring-beam assembly}
\label{subsec:Sensitivityanalysis}

\begin{figure}
\centering
\includegraphics[width=0.5\textwidth]{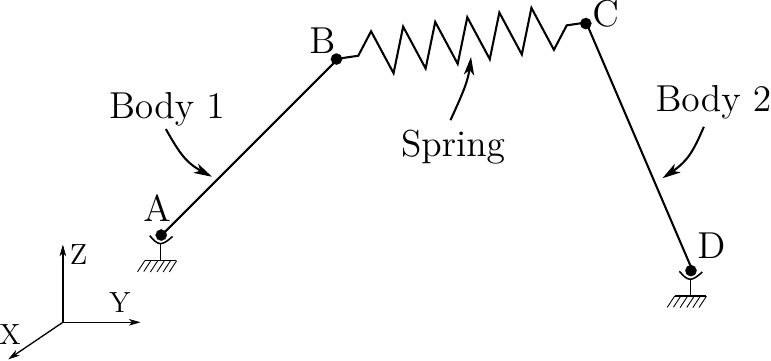}
\caption{\label{fig:case1}A rigid-spring-beam assembly}
\end{figure}

To demonstrate the techniques developed so far, the sensitivity analysis of a simple rigid-spring-beam multibody is presented. The same computations can be employed for more complex assemblies. Suppose the aim is to minimize the squared norm of the displacement of Node B for the mechanism depicted in Figure \ref{fig:case1}.  The objective function in this case is
\begin{equation}\label{eq:objCase1}
\begin{aligned}
\phi=\int_0^T \left\| \vec{\mathbf{d}}_{\mathrm{B}} \right\|^2 dt
\end{aligned} 
\end{equation}
where $\vec{\mathbf{d}}_{\mathrm{B}}$ represents the time-dependent displacement of Node B. The vector of design variables contains three geometrical variables: the initial global position of Node D along X, Y and Z axes. Body 1 is a rigid bar and Body 2 is a beam discretized by one element with length $l$ and constant cross section $A$. The assembly is connected to the ground via two spherical joints at Points A and D. If $X_{\mathrm{B}}$, $Y_{\mathrm{B}}$ and $Z_{\mathrm{B}}$ denote the initial global position of Node B at $t=0$, $\vec{\mathbf{d}}_{\mathrm{B}}$ reads as

\begin{equation}\label{eq:disCase1}
\begin{aligned}
\vec{\mathbf{d}}_{\mathrm{B}} = \mathbf{S}_{\mathrm{Body 1}}\left. \right|_{\mathrm{B}}\vec{\mathbf{q}}_{\mathrm{Body 1}}(t)-\begin{Bmatrix}
X_{\mathrm{B}} \\ Y_{\mathrm{B}} \\ Z_{\mathrm{B}}
\end{Bmatrix}
\end{aligned} 
\end{equation}

Following Equation \ref{eq:disObj}, Equation \ref{eq:objCase1} can be transformed to

\begin{equation}\label{eq:disObjCase1}
\Phi \left( \vec{\mathbf{q}}_{0_{\mathrm{Body 1}}}, \vec{\mathbf{q}}_{1_{\mathrm{Body 1}}}, \dots, \vec{\mathbf{q}}_{N_{\mathrm{Body 1}}} \right) = \displaystyle\sum_{n=0}^{N-1} \left(h \left\|\vec{\mathbf{d}}_{\mathrm{B}} \left( (1-\alpha)\vec{\mathbf{q}}_{n_{\mathrm{Body 1}}} + \alpha\vec{\mathbf{q}}_{{(n+1)}_{\mathrm{Body 1}}} \right) \right\|^2 \right)
\end{equation}

Based on the notation used in Equation \ref{eq:motionEq} and pursuing the equations derived in Sections \ref{sec:Dynamicsof}-\ref{sec:SensitivityGeo}, with $\alpha=0.5$, the set of motion equations for this assembly is
	
\begin{equation}\label{eq:c0cn}
\begin{aligned}
\vec{\mathbf{c}}_0  &:= \frac{1}{h} \mathbf{M} \left( \vec{\mathbf{q}}_1 - \vec{\mathbf{q}}_0 \right) + \frac{h}{2} \frac{\partial}{\partial \vec{\mathbf{q}}_0} U\left(\frac{\vec{\mathbf{q}}_0+ \vec{\mathbf{q}}_1}{2}\right) + \frac{h}{2} \left[ \frac{\partial \vec{\mathbf{g}} (\vec{\mathbf{q}}_0)}{\partial \vec{\mathbf{q}}_0} \right]^{\mathrm{T}} \vec{\boldsymbol{\uplambda}}_0 - \mathbf{M} \dot{\vec{\mathbf{q}}}_0 = \vec{\mathbf{0}}\\[0.5em]
\vec{\mathbf{c}}_n &:= \frac{1}{h} \mathbf{M} \left( 2\vec{\mathbf{q}}_n - \vec{\mathbf{q}}_{n+1} - \vec{\mathbf{q}}_{n-1} \right) - \frac{h}{2} \left( \frac{\partial}{\partial \vec{\mathbf{q}}_n} U\left(\frac{\vec{\mathbf{q}}_n+ \vec{\mathbf{q}}_{n+1}}{2} \right) + \frac{\partial}{\partial \vec{\mathbf{q}}_n} U\left(\frac{\vec{\mathbf{q}}_{n-1}+ \vec{\mathbf{q}}_{n}}{2} \right) \right) \\
&- \frac{h}{2} \left[ \frac{\partial \vec{\mathbf{g}} (\vec{\mathbf{q}}_n)}{\partial \vec{\mathbf{q}}_n} \right]^{\mathrm{T}} \vec{\boldsymbol{\uplambda}}_n = \vec{\mathbf{0}} \quad \text{for $n=1, 2, \dotsc, N-1$}\\
\end{aligned}
\end{equation}

where

\begin{equation}\label{eq:constU}
\begin{gathered}
\vec{\mathbf{q}}_n=\begin{Bmatrix} \vec{\mathbf{q}}_{n_{\mathrm{Body1}}} \\[0.5em] \vec{\mathbf{q}}_{n_{\mathrm{Body2}}} \end{Bmatrix}, \quad \mathbf{M} = \begin{bmatrix} \mathbf{M}_{\text{Body1}} & \mathbf{0} \\[1.0em]
\mathbf{0} & \mathbf{M}_{\text{Body2}} \end{bmatrix} \\
\vec{\mathbf{g}}_{n} = \vec{\mathbf{g}}(\vec{\mathbf{q}}_n) = \begin{bmatrix} \mathbf{S}_{\text{Body1}} \Big|_{\mathrm{A}} & \mathbf{0} \\[1.0em]
\mathbf{0} & \mathbf{S}_{\text{Body2}} \Big|_{\mathrm{D}} \end{bmatrix} \vec{\mathbf{q}}_n - \begin{Bmatrix}
X_{\mathrm{A}} \\ Y_{\mathrm{A}} \\ Z_{\mathrm{A}} \\ X_{\mathrm{D}} \\ Y_{\mathrm{D}} \\ Z_{\mathrm{D}} \end{Bmatrix} = \vec{\mathbf{0}} \\[0.5em]
\frac{\partial U}{\partial \vec{\mathbf{q}}_n} = \frac{\partial U_{\text{gravity}}}{\partial \vec{\mathbf{q}}_n} + \frac{\partial U_{\text{spring}}}{\partial \vec{\mathbf{q}}_n} + \frac{\partial U_{\text{elastic}}}{\partial \vec{\mathbf{q}}_n}
\end{gathered}
\end{equation}

These equations, along with the adjoint equations (Equations \ref{eq:adjointsystem1} and \ref{eq:adjointsystem2}), as well as Equations \ref{eq:position}-\ref{eq:dampereq} can be utilized to find adjoint variables $\vec{\boldsymbol{\upmu}}_n$ and $\vec{\boldsymbol{\upeta}}_n$. To compute the total derivative based on Equation \ref{eq:finalSensitivity}, $\frac{\partial \Phi}{\partial a_i}$, $\frac{\partial \vec{\mathbf{c}}_n}{\partial a_i}$ and $\frac{\partial \vec{\mathbf{g}}_n}{\partial a_i}$ are also required. If $X_{\mathrm{D}}$, $Y_{\mathrm{D}}$ and $Z_{\mathrm{D}}$ are the initial global position of Node D at $t=0$, since there is no explicit dependency in $\Phi$ on the given design variables, $\frac{\partial \Phi}{\partial X_{\mathrm{D}}}=\frac{\partial \Phi}{\partial Y_{\mathrm{D}}}=\frac{\partial \Phi}{\partial Z_{\mathrm{D}}}=0$. As for the constraint equations, there are two spherical joints in this example, only one of which is a function of $X_{\mathrm{D}}$, $Y_{\mathrm{D}}$ and $Z_{\mathrm{D}}$. Using Equation \ref{eq:constU}

\begin{equation}\label{eq:DjointCase1}
\begin{aligned}
\frac{\partial \vec{\mathbf{g}}_n}{\partial X_{\mathrm{D}}} = - \begin{Bmatrix}
0 \\ 0 \\ 0 \\ 1 \\ 0 \\ 0 \end{Bmatrix}, \quad \frac{\partial \vec{\mathbf{g}}_n}{\partial Y_{\mathrm{D}}} = - \begin{Bmatrix}
0 \\ 0 \\ 0 \\ 0 \\ 1 \\ 0 \end{Bmatrix}, \quad \frac{\partial \vec{\mathbf{g}}_n}{\partial Z_{\mathrm{D}}} = - \begin{Bmatrix}
0 \\ 0 \\ 0 \\ 0 \\ 0 \\ 1 \end{Bmatrix}
\end{aligned}
\end{equation}

In the equations of $\vec{\mathbf{c}}_0$ and $\vec{\mathbf{c}}_n$, Equation \ref{eq:c0cn}, only the terms related to Body 2 are explicitly dependent on $X_{\mathrm{D}}$, $Y_{\mathrm{D}}$ and $Z_{\mathrm{D}}$. They are

\begin{itemize}
  \item $\mathbf{M}_{\mathrm{Body2}}$
  
  According to Equation \ref{eq:kineticdeltaU}
  \begin{equation}\label{eq:Mbeam}
  \begin{aligned}
  \frac{\partial \mathbf{M}_{\mathrm{Body2}}}{\partial a_i} &= \frac{\partial \mathbf{M}_{\mathrm{Body2}}}{\partial l} \frac{\partial l}{\partial a_i}
  \end{aligned}
  \end{equation}
  \item $\partial U_{\text{gravity}_\mathrm{Body2}} / \partial \vec{\mathbf{q}}_n$
  
  Denoting $\vec{\boldsymbol{\upnu}}$ as the gravity vector and following Equation \ref{eq:gengravbeam}
  \begin{equation}
  \begin{aligned}
  \frac{\partial}{\partial a_i} \left[\frac{\partial U_{\text{gravity}_\mathrm{Body2}}}{\partial \vec{\mathbf{q}}_n} \right] &= -\frac{\partial}{\partial l} \left[\int\limits_V \left[\rho \vec{\boldsymbol{\upnu}}^{\mathrm{T}}\mathbf{S} \right]dV \right] \frac{\partial l}{\partial a_i} 
  \end{aligned}
  \end{equation}
  \item $\partial U_{\text{longitudinal}_\mathrm{Body2}} / \partial \vec{\mathbf{q}}_n$ and $\partial U_{\text{transverse}_\mathrm{Body2}} / \partial \vec{\mathbf{q}}_n$
  
  Assuming small longitudinal deformation and utilizing Equation \ref{eq:linearForm}
  \begin{equation}\label{eq:DUBeam}
  \begin{aligned}
  \frac{\partial}{\partial a_i} \left[\frac{\partial U_{\text{longitudinal}_\mathrm{Body2}}}{\partial \vec{\mathbf{q}}_n} \right] &= \frac{\partial \mathbf{K}_{\text{longitudinal}}}{\partial l} \frac{\partial l}{\partial a_i} \\
  \frac{\partial}{\partial a_i} \left[\frac{\partial U_{\text{longitudinal}_\mathrm{Body2}}}{\partial \vec{\mathbf{q}}_n} \right] &= \frac{\partial \mathbf{K}_{\text{transverse}}}{\partial l} \frac{\partial l}{\partial a_i}
  \end{aligned}
  \end{equation} 
\end{itemize}

In Equations \ref{eq:Mbeam}-\ref{eq:DUBeam}, $l$ is the beam's length in the undeformed configuration, which is calculated by

\begin{equation}
l=\sqrt{\left(X_{\mathrm{D}}-X_{\mathrm{C}} \right)^2+\left(Y_{\mathrm{D}}-Y_{\mathrm{C}} \right)^2+\left(Z_{\mathrm{D}}-Z_{\mathrm{C}} \right)^2}
\end{equation}

and its derivatives are

\begin{equation}
\frac{\partial l}{\partial X_{\mathrm{D}}} = \frac{X_{\mathrm{D}}-X_{\mathrm{C}}}{l}, \quad \frac{\partial l}{\partial Y_{\mathrm{D}}} = \frac{Y_{\mathrm{D}}-Y_{\mathrm{C}}}{l}, \quad \frac{\partial l}{\partial Z_{\mathrm{D}}} = \frac{Z_{\mathrm{D}}-Z_{\mathrm{C}}}{l}
\end{equation}

Having computed all these terms and substituting them into Equation \ref{eq:finalSensitivity}, the sensitivity of the given objective function is thus  found. For cases with a higher number of components and other types of geometrical design variables, derivatives are obtainable in a similar manner.

\section{Numerical Examples}
\label{sec:numerical}

{Three numerical test cases are provided in this section. The first example is to validate the proposed sensitivity analysis technique based on the study reported in{\cite{Pi2012}}. The other two are to compare its performance with the finite difference method (FD) and assay its applicability to the cases containing geometrical and non-geometrical design variables. In all examples, a Newton-Raphson scheme is used to solve the linear and nonlinear equations of forward and backward simulations. Once the required gradients are computed, the NLopt library{\cite{johnson2014nlopt}}, which is an open-source nonlinear optimization library, is utilized to find the optimum solutions. To handle the optimization constraints, the Augmented Lagrangian method is exploited. In all examples, $\alpha$ is set to 0.5 in the forward and backward simulations. The optimization scheme is stopped when the objective function improvement is less than $10^{-6}$ for five consecutive iterations or the maximum number of iterations, set to 60, is reached.}

\subsection{Sensitivity analysis of a flexible pendulum}
\label{subsec:Case1}
\begin{figure}[H]
\centering
{\includegraphics[width=0.6\textwidth]{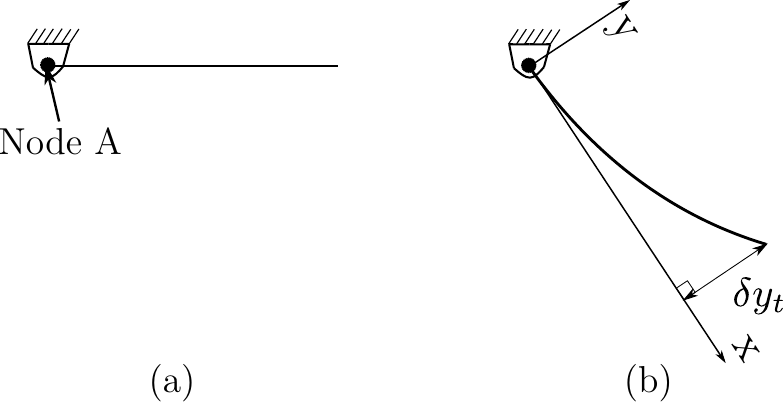}}
\caption{\label{fig:pendulum}{(a) Pendulum at $t=0$ (b) Deflected pendulum at time instance $t$}}
\end{figure}

{Consider the planar flexible pendulum shown in Figure {\ref{fig:pendulum}}. It is made of a material with the density $\rho = 4000 \: \mathrm{kg/m^3}$, Young's modulus $E=10^7 \: \mathrm{N/m^2}$ and Poisson's ratio 0.3. The length of the beam is 1.2 m, discretized by 5 equi-length elements, and the width ($h$) of its square cross-section is $0.05 \: \mathrm{m}$. The duration of the problem is 4 seconds and the beam is initially at rest horizontally along the global X axes.}

\begin{figure}[H]
\centering
{\includegraphics[width=0.9\textwidth]{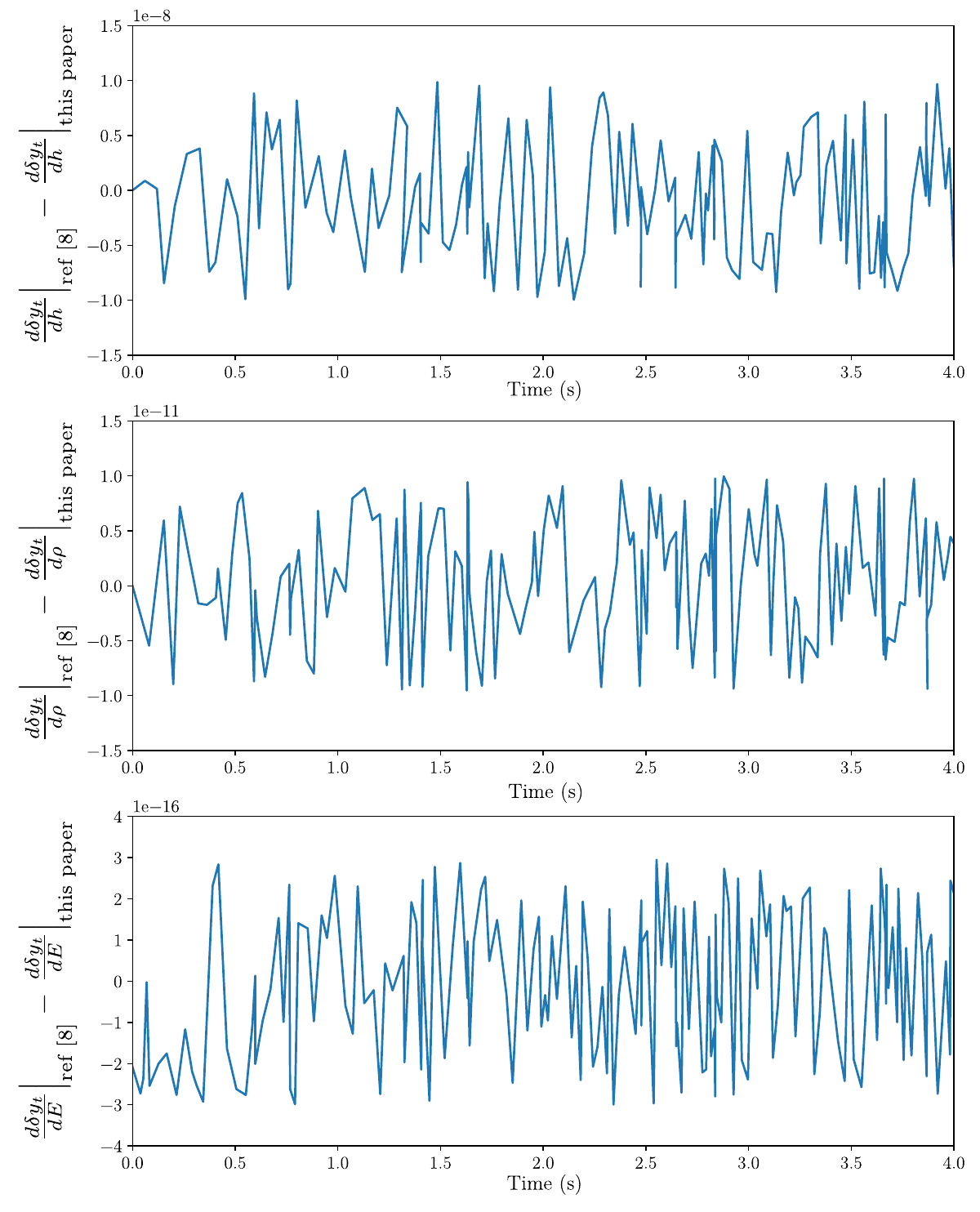}}
\caption{\label{fig:pendulumsen} {Difference between the sensitivity values with respect to $h$, $\rho$ and $E$ from{\cite{Pi2012}} and the proposed DAVM}}
\end{figure}

{Suppose $\delta y_t$ is the deflection of the tip at time $t$, and the goal is to compute the sensitivity of this parameter at different time instances with respect to $h$, $\rho$ and $E$. $\delta y_t$ is measured along the y axis of the local coordinate system whose origin is at Node A and its x axis is tangent to the beam's centerline. Figure {\ref{fig:pendulumsen}} shows the difference between the values reported in{\cite{Pi2012}} and those computed using the developed DAVM in this paper. The results are almost identical, and the small differences are due to the different time-stepping methods for solving the forward and backward simulation equations.}

\subsection{Optimization of a rigid-spring-beam assembly}
\label{subsec:Case2}

For the second example, the two-dimensional version of the assembly discussed in Section \ref{sec:SensitivityGeo} and shown in Figure \ref{fig:case1} is considered. The system is connected to the ground via two spherical joints at Nodes A and D. Body 1 is a rigid cylinder connecting Nodes A and B. Its mass and cross-sectional diameter are 1.6 kg and 0.01 m, respectively, and its center of mass is located in the midpoint between the two ends. Body 2 is a flexible beam with Young's modulus of 70 GPa, Poisson's ratio of 0.3 and density of 3700 $\mathrm{kg/m^3}$. Its cross-section is a solid square with 0.01 m width. A linear spring connects Nodes B and C. The spring constant is 100 $\mathrm{N/m^3}$ and its initial undeformed length is equal to the distance between Nodes B and C at time zero. 

In this example, the duration is 1 sec, and the problem is run by two different time-step sizes of 0.001 sec and 0.0005 sec. To drive the mechanism an initial angular velocity of -5 rad/sec about the positive Z axis is applied to Body 2. The objective function is the squared norm of Node B's displacement over the entire simulation duration. A lower bound constraint is imposed on the beam's length to prevent producing a beam with zero length. This optimization problem can be formulated as

\begin{equation}\label{eq:obj1}
\begin{gathered}
\min_{X_{\mathrm{D}}, Y_{\mathrm{D}}} \int_0^1 \left\| \vec{\mathbf{d}}_{\mathrm{B}} \right\|^2 dt \\
\hspace{-18em}\text{subject to} \\[0.5em]
l_{\text{Body 2}} \geq 0.001 \mathrm{\:m}, \\[0.5em]
-10 \mathrm{\:m} \leq X_{\mathrm{D}}, Y_{\mathrm{D}} \leq 10 \mathrm{\:m}
\end{gathered}
\end{equation}

\begin{figure}[H]
\centering
\includegraphics[width=0.8\textwidth]{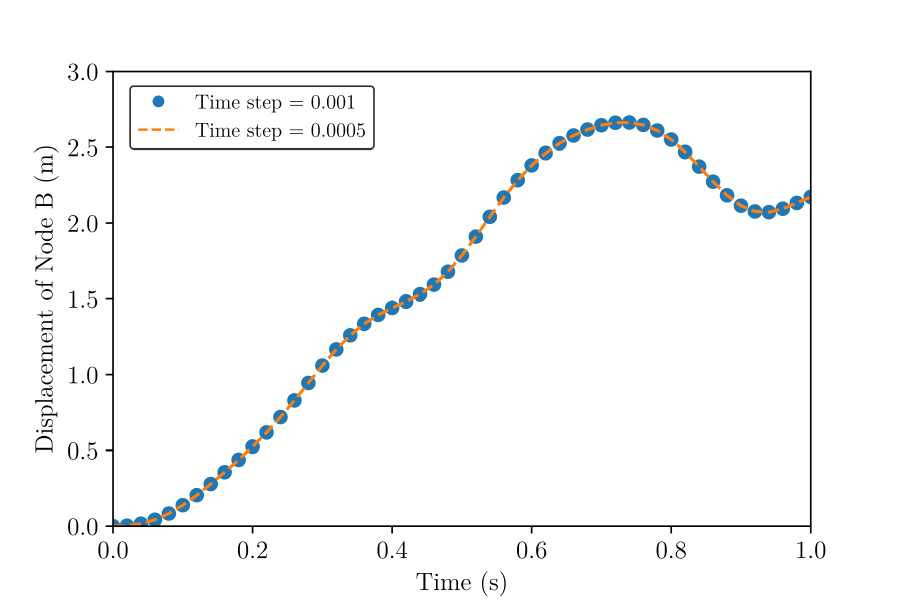}
\caption{\label{fig:TwoDis}Displacement of Node B for the two time-steps}
\end{figure}

Figure \ref{fig:TwoDis} shows the displacement of Node B in time for the initial configuration for the two given time-step sizes. As can be seen, the plots for the two time-steps coincide and so the forward simulation phase is insensitive to these two time-step sizes.

\begin{figure}[H]
\centering
\includegraphics[width=0.8\textwidth]{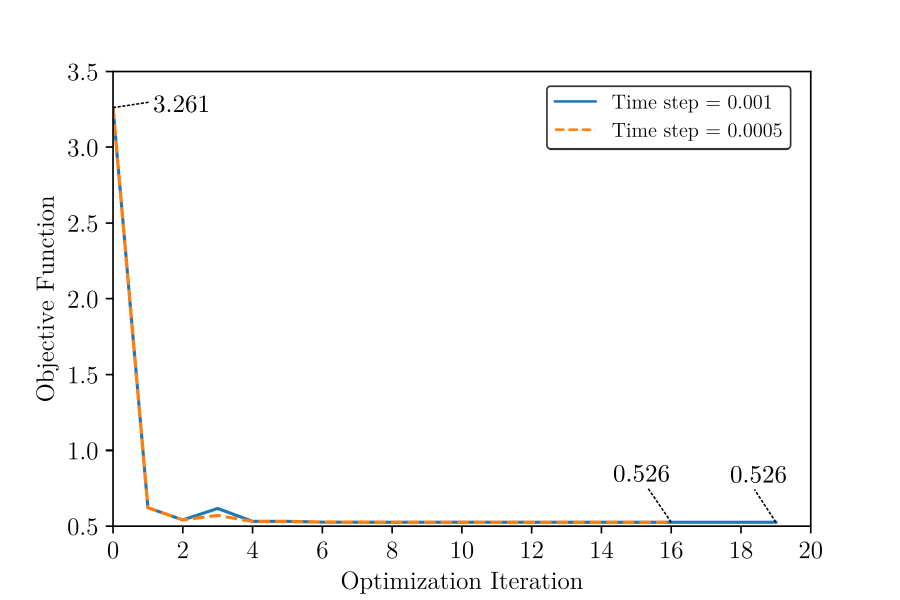}
\caption{\label{fig:ObjEvo}Evolution of objective function for the two time-steps}
\end{figure}

\begin{figure}[H]
\centering
\includegraphics[width=0.8\textwidth]{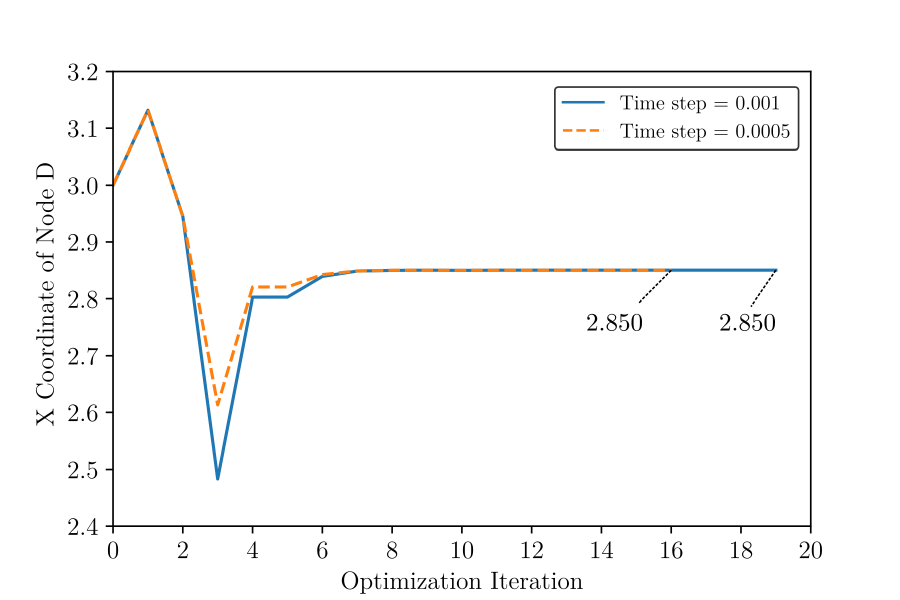}
\caption{\label{fig:XEvo}Evolution of X coordinate of Node D for the two time-steps}
\end{figure}

\begin{figure}[H]
\centering
\includegraphics[width=0.8\textwidth]{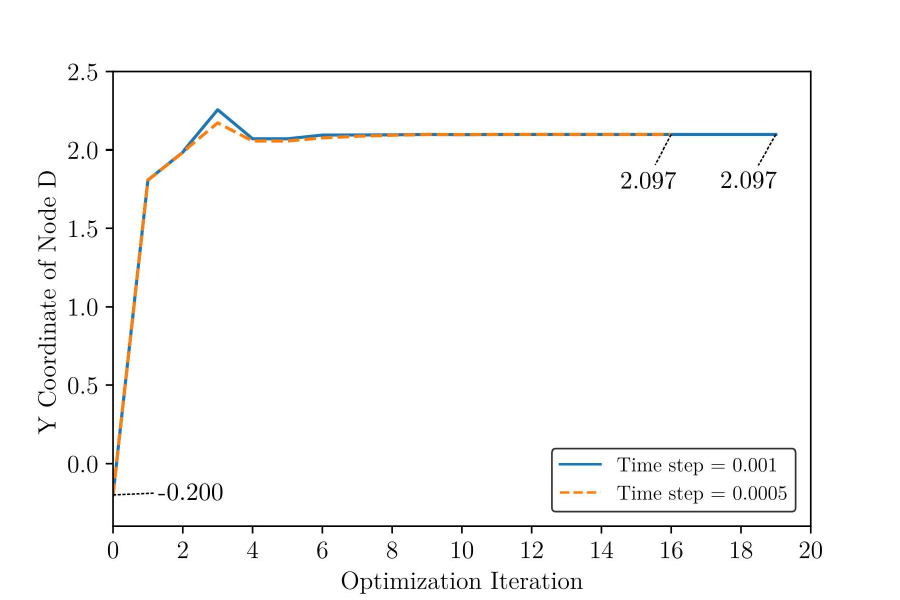}
\caption{\label{fig:YEvo}Evolution of Y coordinate of Node D for the two time-steps}
\end{figure}

{Figures {\ref{fig:ObjEvo}-\ref{fig:YEvo}} present the evolution of the objective function and design variables for both time-step sizes. Even though the optimization procedure proceeds slightly differently for the two time-steps, the final optimum solutions are almost exactly identical. These figures confirm that not only the forward simulation phase, but also the backward simulation and consequently the sensitivity values and optimization results are independent of the time-step size.  This is, however, not the case in the approximation methods, such as the equivalent static loads (ESL). In ESL, the sensitivity values are computed for a sequence of equivalent linear static models which have the same response fields as those in the actual dynamic problem at a selected number of time instances. Depending on which time-steps are taken into account and for how many of them the equivalent model is constructed, the optimum solutions could be entirely different (e.g., the numerical examples in{\cite{Sun2016}}). It has been even shown that ESL does not in general lead to optimal designs{\cite{stolpe2018equivalent}}, whereas the proposed sensitivity analysis technique as explained in Section {\ref{sec:problemdef}} guarantees moving along the descent direction (or ascent direction in maximization problems) and subsequently improving the defined objective function.}

\begin{figure}[H]
\centering
\includegraphics[width=0.8\textwidth]{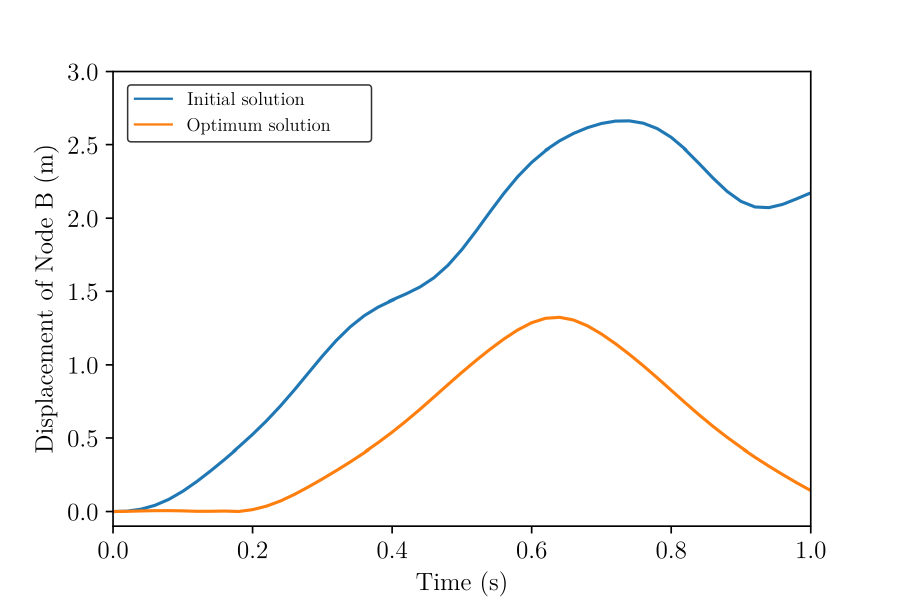}
\caption{\label{fig:OptNotOpt}Displacement of Node B in the initial and optimum configurations}
\end{figure}

Figure \ref{fig:OptNotOpt} shows the displacement of Node B in the initial and optimum configurations. According to this figure, the displacement of Node B in the optimum solution is considerably smaller than that in the initial solution for all time steps, therefore leading to a significantly better performance in terms of the defined criterion in Equation \ref{eq:obj1}. 

\begin{figure}[H]
\centering
{\includegraphics[width=0.8\textwidth]{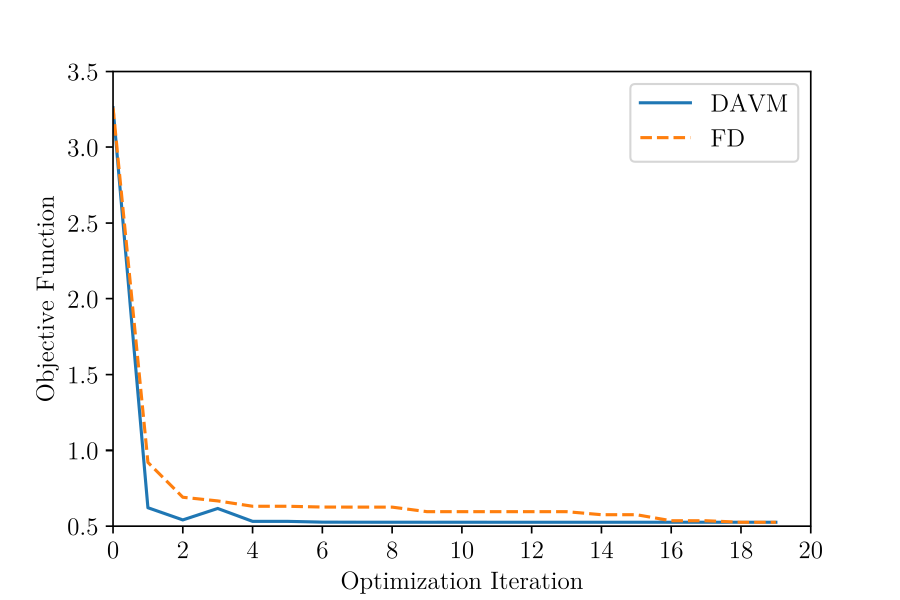}}
\caption{\label{fig:DAVMFD} {Evolution of objective function using DAVM and FD for $\Delta t = 0.001$ sec}}
\end{figure}

{To further investigate the performance of the proposed sensitivity analysis technique, this problem is run using FD as well. Figure {\ref{fig:DAVMFD}} depicts the optimization evolution for the developed DAVM and FD with 1 percent perturbation value. Even though both approaches converge to the same optimum solution, their computation time is considerably different. On the same computer, the DAVM executes about 4.5 times faster than FD. This is mainly due to the fact that in FD to compute the sensitivities for each design variable, an additional round of simulation is required. The difference is even more significant for more complex examples, which makes FD almost impractical in those cases.}

\newpage
\subsection{Optimization of an automotive double-wishbone suspension system}
\label{subsec:Case3}

\begin{figure}[!htbp]
\centering
\includegraphics[width=0.55\textwidth]{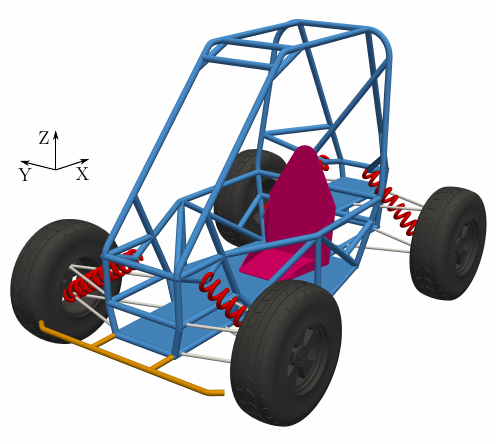}
\caption{\label{fig:SAE}The car with double-wishbone suspension system}
\end{figure}

\begin{figure}
\centering
\includegraphics[width=0.6\textwidth]{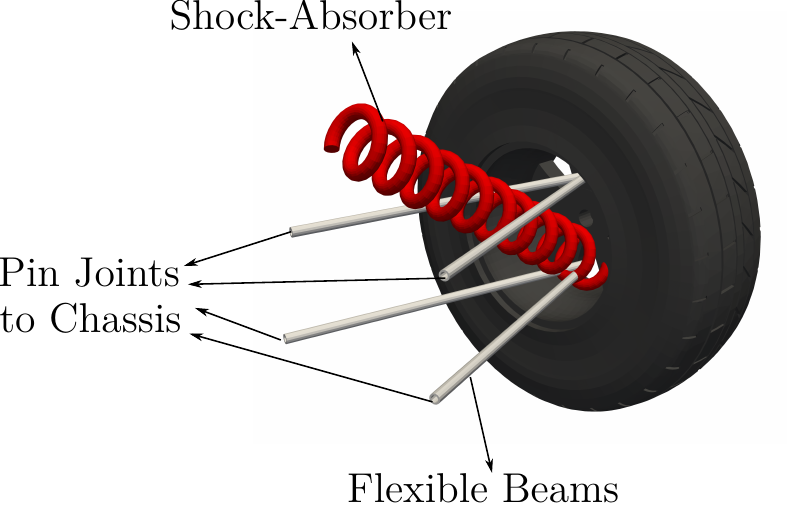}
\caption{\label{fig:susModule}A suspension system module}
\end{figure}

To examine the developed methods on a large-scale engineering case study, they are applied to optimize the double-wishbone suspension system of the car depicted in Figure \ref{fig:SAE}. There are four suspension system modules, as portrayed in Figure \ref{fig:susModule}, each consists of two wish-bone like structures and a shock-absorber.  The wish-bones are composed of two hollow circular flexible beams which are welded together from one end and connected to the wheels from another end through spherical joints. The shock-absorbers are modeled by a spring and a damper, connected to the wheels and the chassis. The rest of the bodies are considered as rigid bodies. No symmetry is imposed on the design in order to give the optimization scheme the freedom to find the optimum values of design variables independently for each of four suspension system modules. 

\begin{figure}
\centering
\includegraphics[width=0.8\textwidth]{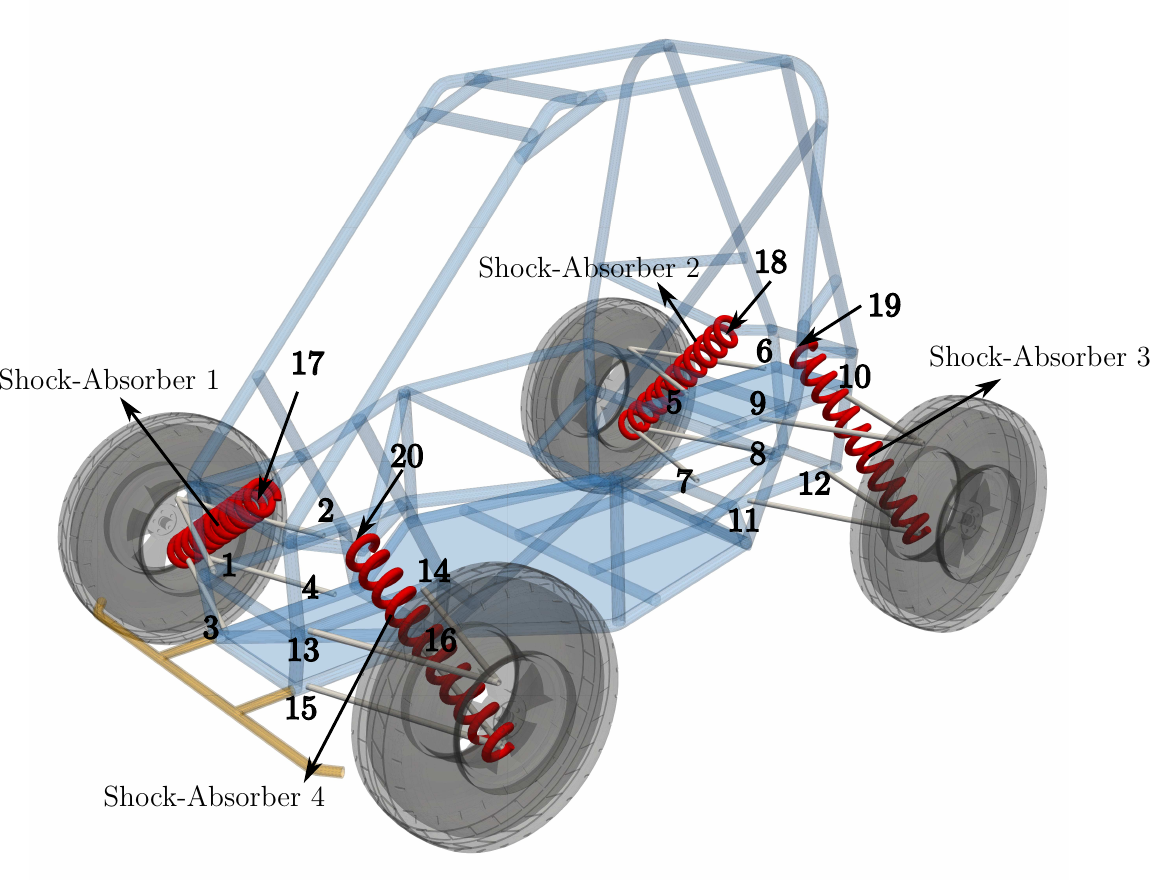}
\caption{\label{fig:SAETrans}Node IDs associated with design variables}
\end{figure}

The objective function is the squared magnitude of displacement of the chassis’s center of mass over time. There are 52 geometrical and 8 non-geometrical design variables in this problem. The geometrical variables are the position of the connection points between the beams and chassis ($X$, $Y$ and $Z$ coordinates of each point) and those between the shock-absorbers and chassis along the global Z axis (see Figure \ref{fig:SAETrans}). The non-geometrical variables are spring constants and damping coefficients of each shock-absorber. The lower and upper bounds and initial values of each parameter are reported in Table \ref{table:dvar}. 

\begin{table}
\caption{\label{table:dvar}Initial values and range of design variables for suspension system example}
\begin{center}
\resizebox{\textwidth}{!}{%
\begin{tabular}{llll}
\toprule
\\[-0.7em]
\textbf{Design variables ($a_i$)}   & \textbf{Initial values}  & \textbf{Lower bounds ($lb_i$)}  & \textbf{Upper bounds ($ub_i$)}  \\[0.5em]
\midrule
$\left(X_1,Y_1,Z_1 \right)$      & $(-1.270, 0.210, 0.150)$ m    & $(-1.350, 0.100, 0.100)$  & $(-1.100, 0.350, 0.300)$     \\[0.5em]
$\left(X_2,Y_2,Z_2 \right)$      & $(-0.960, 0.210, 0.150)$ m    & $(-1.100, 0.100, 0.100)$  & $(-0.800, 0.350, 0.300)$     \\[0.5em]
$\left(X_3,Y_3,Z_3 \right)$      & $(-1.270, 0.180, 0)$ m    & $(-1.350, 0.100, -0.100)$  & $(-1.100, 0.350, 0.100)$     \\[0.5em]
$\left(X_4,Y_4,Z_4 \right)$      & $(-0.960, 0.180, 0)$ m    & $(-1.100, 0.100, -0.100)$  & $(-0.800, 0.350, 0.100)$     \\[0.5em]
$\left(X_5,Y_5,Z_5 \right)$      & $(0.210, 0.140, 0.250)$ m    & $(0.100, 0.100, 0.100)$  & $(0.350, 0.300, 0.300)$     \\[0.5em]
$\left(X_6,Y_6,Z_6 \right)$      & $(0.520, 0.140, 0.250)$ m    & $(0.400, 0.100, 0.100)$  & $(0.650, 0.300, 0.300)$     \\[0.5em]
$\left(X_7,Y_7,Z_7 \right)$      & $(0.210, 0.100, 0)$ m    & $(0.100, 0.100, -0.100)$  & $(0.350, 0.300, 0.100)$     \\[0.5em]
$\left(X_8,Y_8,Z_8 \right)$      & $(0.520, 0.100, 0)$ m    & $(0.400, 0.100, -0.100)$  & $(0.650, 0.300, 0.100)$     \\[0.5em]
$\left(X_9,Y_9,Z_9 \right)$      & $(0.210, -0.140, 0.250)$ m    & $(0.100, -0.300, 0.100)$  & $(0.350, -0.100, 0.300)$     \\[0.5em]
$\left(X_{10},Y_{10},Z_{10} \right)$      & $(0.520, -0.140, 0.250)$ m    & $(0.400, -0.300, 0.100)$  & $(0.650, -0.100, 0.300)$     \\[0.5em]
$\left(X_{11},Y_{11},Z_{11} \right)$      & $(0.210, -0.100, 0)$ m    & $(0.100, -0.300, -0.100)$  & $(0.350, -0.100, 0.100)$     \\[0.5em]
$\left(X_{12},Y_{12},Z_{12} \right)$      & $(0.520, -0.100, 0)$ m    & $(0.400, -0.300, -0.100)$  & $(0.650, -0.100, 0.100)$     \\[0.5em]
$\left(X_{13},Y_{13},Z_{13} \right)$      & $(-1.270, -0.210, 0.150)$ m    & $(-1.350, -0.350, 0.100)$  & $(-1.100, -0.100, 0.300)$     \\[0.5em]
$\left(X_{14},Y_{14},Z_{14} \right)$      & $(-0.960, -0.210, 0.150)$ m    & $(-1.100, -0.350, 0.100)$  & $(-0.800, -0.100, 0.300)$     \\[0.5em]
$\left(X_{15},Y_{15},Z_{15} \right)$      & $(-1.270, -0.180, 0)$ m    & $(-1.350, -0.350, -0.100)$  & $(-1.100, -0.100, 0.100)$     \\[0.5em]
$\left(X_{16},Y_{16},Z_{16} \right)$      & $(-0.960, -0.180, 0)$ m    & $(-1.100, -0.350, -0.100)$  & $(-0.800, -0.100, 0.100)$     \\[0.5em]
$Z_{17}$	&	0.3 m	&	0.15	&	0.6		\\[0.5em]
$Z_{18}$	&	0.4 m	&	0.15	&	0.6		\\[0.5em]
$Z_{19}$	&	0.4 m	&	0.15	&	0.6		\\[0.5em]
$Z_{20}$	&	0.3 m	&	0.15	&	0.6		\\[0.5em]
$k_{1}$	&	10000 N/m	&	100		&	30000	\\[0.5em]
$k_{2}$	&	10000 N/m	&	100		&	30000	\\[0.5em]
$k_{3}$	&	10000 N/m	&	100		&	30000	\\[0.5em]
$k_{4}$	&	10000 N/m	&	100		&	30000	\\[0.5em]
$c_{1}$	&	50 N.s/m	&	20		&	1000	\\[0.5em]
$c_{2}$	&	500 N.s/m	&	20		&	1000	\\[0.5em]
$c_{3}$	&	500 N.s/m	&	20		&	1000	\\[0.5em]
$c_{4}$	&	50 N.s/m	&	20		&	1000	\\[0.5em]
\bottomrule
\end{tabular}}
\end{center}
\end{table}

\clearpage
To mimic the driving condition on a rugged terrain, the front and rear wheels are subject to two different harmonic excitations along the global Z axis as described in Table \ref{table:excitation}. For this problem, the final time $T$ and the time-step size are 3 sec and 0.001 sec, respectively. Other parameters are summarized in Table \ref{table:paramCase2}. 

\begin{table}
\caption{\label{table:excitation}Front and rear wheels excitation}
\begin{center}
\begin{tabular}{ll}
\toprule
Front wheels \hspace{5cm}  & $d_z=0.050 \sin \left(2\pi t \right)$  \\[0.5em]
Rear wheels  \hspace{5cm}  & $d_z=0.075 \sin \left(4\pi t \right)$  \\
\bottomrule
\end{tabular}
\end{center}
\end{table}

\begin{table}[h]
\caption{\label{table:paramCase2}Physical and material properties for suspension system example}
\begin{center}
\begin{tabular}{ll}
\toprule
\\[-0.7em]
\textbf{Parameters} \hspace{7cm}  & \textbf{Values}  \\[0.5em]
\midrule
mass of chassis & 300 kg \\[0.5em]
mass of wheels & 20 kg \\[0.5em]
Young’s modulus of beams & 200 GPa \\[0.5em]
Poisson’s ratio of beams & 0.3 \\[0.5em]
outer radius of beams’ cross-section & 0.015 m \\[0.5em]
inner radius of beams’ cross-section & 0.01 m \\[0.5em]
density of beams & 7800 kg/$\mathrm{m}^3$ \\[0.5em]
location of chassis’s center of mass & $(0.243, 0, 0.365)$ m \\[0.5em]
location of front-left wheel’s center of mass & $(-1.110, -0.680, 0.075)$ m \\[0.5em]
location of front-right wheel’s center of mass & $(-1.110, 0.680, 0.075)$ m \\[0.5em]
location of rear-left wheel’s center of mass & $(0.370, -0.560, 0.125)$ m \\[0.5em]
location of rear-right wheel’s center of mass & $(0.370, 0.560, 0.125)$ m \\
\bottomrule
\end{tabular}
\end{center}
\end{table}

\clearpage
The optimization problem in this case is written as

\begin{equation}
\begin{gathered}
\min_{a_i \in \mathbb{R}^{60}} \int_0^3 \left\| \vec{\mathbf{d}}_{CM} \right\|^2 dt \\
\hspace{-18em}\text{subject to} \\[0.5em]
lb_i \leq a_i \leq ub_i \quad (i=1,\dots,60), \\[0.5em]
-400 \mathrm{\:MPa} \leq \sigma_{beam_j} \leq 400 \mathrm{\:MPa} \quad (j=1,\dots,16), \\[0.5em]
l_{beam_j} \geq 0.001 \mathrm{\:m} \quad (j=1,\dots,16)
\end{gathered}
\end{equation}

\begin{figure}[!htbp]
\centering
{\includegraphics[width=0.8\textwidth]{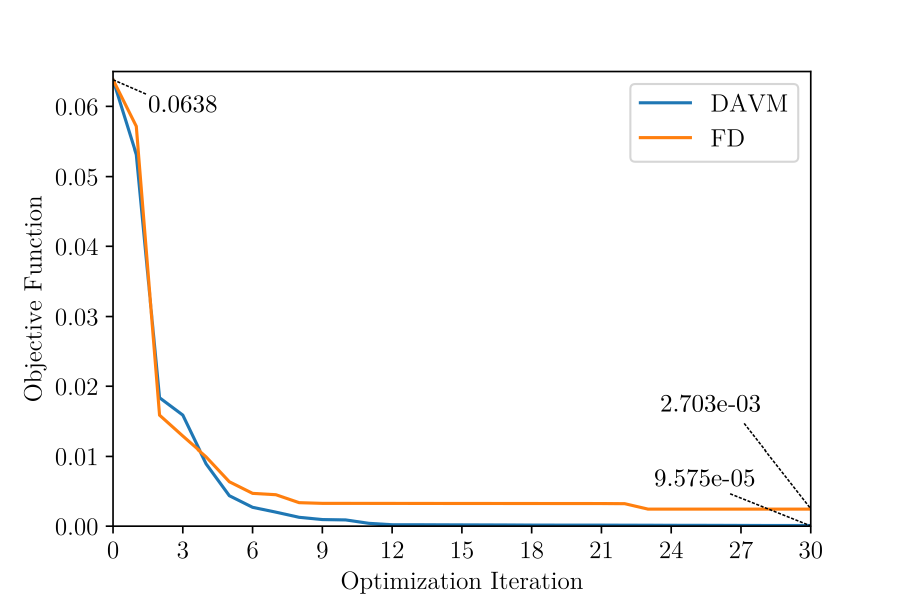}}
\caption{\label{fig:ObjSus} {Evolution of objective function for the suspension system example}}
\end{figure}

{Similar to the previous example, for the sensitivity analysis, FD with 1 percent perturbation value is employed as well. Figure {\ref{fig:ObjSus}} shows the evolution of the objective function using the proposed DAVM and FD. Both methods converge to an optimum design, however, the one found by the DAVM has a better objective function value. This could be due to approximating the sensitivity values in FD, as opposed to having the exact values using the DAVM. Errors in computing the gradients drive the optimization procedure along either a wrong direction (i.e., deteriorating the performance) or the one that does not improve the performance the most (i.e., slowing down the convergence). In this example, both methods eventually lead to the same solution, if they are run for 60 iterations. In terms of the computation time, on the same computer, the DAVM runs about 40 times faster than FD for 60 optimization iterations.}

\begin{figure}[H]
\centering
\includegraphics[width=0.8\textwidth]{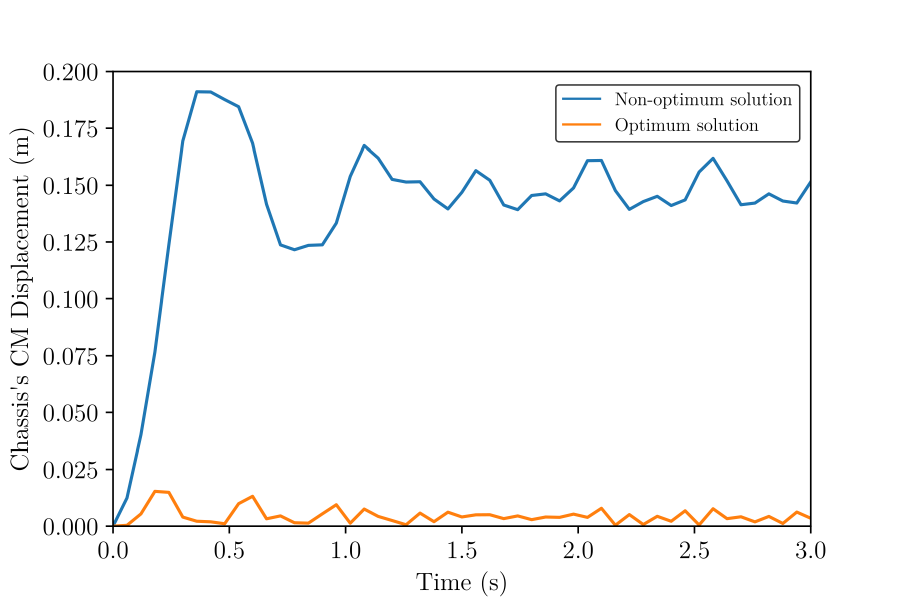}
\caption{\label{fig:OptNotOptSus}Displacement of chassis's center of mass in the non-optimum (initial) and optimum configurations}
\end{figure}

Figure \ref{fig:OptNotOptSus} is the plot of total displacement of the chassis’s center of mass for the initial and optimum designs utilizing the DAVM. It well confirms that the optimum solution outperforms the initial design in terms of the defined performance criterion.

\begin{figure}[H]
\centering
{\includegraphics[width=1.0\textwidth]{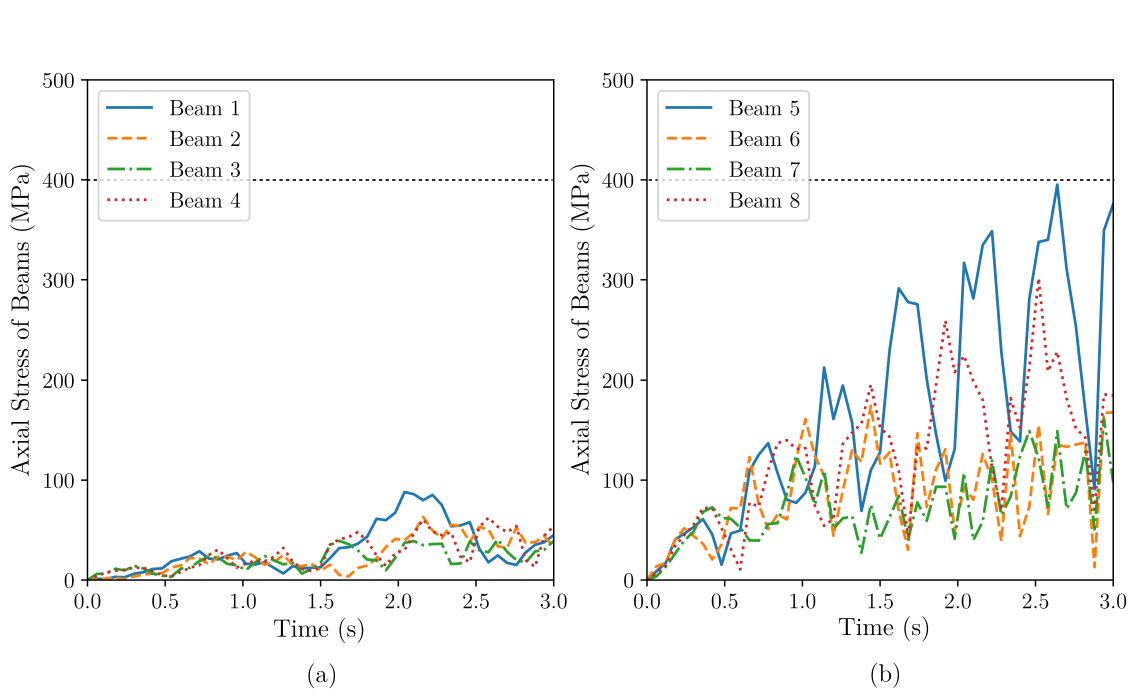}}
\caption{\label{fig:BeamStressSuspension} {Maximum axial stress of beams in the optimum solution (a) Front beams (b) Rear beams}}
\end{figure}

Figure \ref{fig:BeamStressSuspension} shows the variation of beams’ maximum axial stress for the optimum solution. Since the model is symmetric with respect to XZ plane, only the values for one half of the car are plotted. It is clear that all stress constraints are satisfied. Also, as the center of mass of the chassis is closer to the rear of the car, beams belonging to the rear suspension system module (Beams 5, 6, 7, 8) are subject to larger forces and subsequently have higher stress values. 

Optimum values of objective function and design variables using the proposed DAVM and FD are reported in Table \ref{table:Optvar}. Although no symmetry constraint is imposed on the front and rear suspension modules, the optimum values suggest a symmetric solution for each pair, as expected.

\begin{table}
\caption{\label{table:Optvar}Optimum values of design variables for suspension system example}
\begin{center} 
\begin{tabular}{lll}
\toprule
\\[-0.7em]
\textbf{Design variables}    & \textbf{Optimum values(DAVM)} & \textbf{Optimum values(FD)}  \\[0.5em]
\midrule
$\left(X_1,Y_1,Z_1 \right)$      & $(-1.27028, 0.20972, 0.14977)$ m & $(-1.27051, 0.20937, 0.14930)$ m\\[0.5em]
$\left(X_2,Y_2,Z_2 \right)$      & $(-0.96033, 0.21016, 0.15049)$ m & $(-0.96044, 0.21009, 0.15107)$ m\\[0.5em]
$\left(X_3,Y_3,Z_3 \right)$      & $(-1.27009, 0.18006, 0.00561)$ m & $(-1.27005, 0.17924, 0.01282)$ m\\[0.5em]
$\left(X_4,Y_4,Z_4 \right)$      & $(-0.96021, 0.18013, 0.00234)$ m & $(-0.96011, 0.17940, 0.00849)$ m\\[0.5em]
$\left(X_5,Y_5,Z_5 \right)$      & $(0.21030, 0.13594, 0.26894)$ m & $(0.21081, 0.13529, 0.27875)$ m\\[0.5em]
$\left(X_6,Y_6,Z_6 \right)$      & $(0.51827, 0.14297, 0.23080)$ m & $(0.51602, 0.14415, 0.21930)$ m\\[0.5em]
$\left(X_7,Y_7,Z_7 \right)$      & $(0.20898, 0.11675, -0.03229)$ m & $(0.20936, 0.11558, -0.02645)$ m\\[0.5em]
$\left(X_8,Y_8,Z_8 \right)$      & $(0.51810, 0.12266, -0.04791)$ m & $(0.51797, 0.12282, -0.04917)$ m\\[0.5em]
$\left(X_9,Y_9,Z_9 \right)$      & $(0.210301, -0.13596, 0.26894)$ m & $(0.21081, -0.13529, 0.27874)$ m\\[0.5em]
$\left(X_{10},Y_{10},Z_{10} \right)$      & $(0.51827, -0.14299, 0.23080)$ m & $(0.51602, -0.14415, 0.21930)$ m\\[0.5em]
$\left(X_{11},Y_{11},Z_{11} \right)$      & $(0.20898, -0.11675, -0.03229)$ m & $(0.20935, -0.11558, -0.02644)$ m\\[0.5em]
$\left(X_{12},Y_{12},Z_{12} \right)$      & $(0.51810, -0.12265, -0.04791)$ m & $(0.51797, -0.12282, -0.04917)$ m\\[0.5em]
$\left(X_{13},Y_{13},Z_{13} \right)$      & $(-1.27028, -0.20972, 0.14977)$ m & $(-1.27051, -0.20938, 0.14930)$ m\\[0.5em]
$\left(X_{14},Y_{14},Z_{14} \right)$      & $(-0.96033, -0.21016, 0.15049)$ m & $(-0.96044, -0.21009, 0.15107)$ m\\[0.5em]
$\left(X_{15},Y_{15},Z_{15} \right)$      & $(-1.27009, -0.18006, 0.00561)$ m & $(-1.27005, -0.17924, 0.01282)$ m\\[0.5em]
$\left(X_{16},Y_{16},Z_{16} \right)$      & $(-0.96021, -0.18014, 0.00234)$ m & $(-0.96011, -0.17941, 0.00850)$ m\\[0.5em]
$Z_{17}$	&	0.18536 m &	0.1500 m\\[0.5em]
$Z_{18}$	&	0.28891 m &	0.27918 m\\[0.5em]
$Z_{19}$	&	0.28891 m &	0.27918 m\\[0.5em]
$Z_{20}$	&	0.18536 m &	0.1500 m\\[0.5em]
$k_{1}$	&	7552.080 N/m &	3780.350 N/m\\[0.5em]
$k_{2}$	&	21048.100 N/m &	21058.600 N/m\\[0.5em]
$k_{3}$	&	21048.100 N/m &	21058.600 N/m\\[0.5em]
$k_{4}$	&	7552.080 N/m &	3780.350 N/m\\[0.5em]
$c_{1}$	&	102.625 N.s/m &	122.734 N.s/m\\[0.5em]
$c_{2}$	&	457.101 N.s/m &	364.194 N.s/m\\[0.5em]
$c_{3}$	&	457.101 N.s/m &	364.194 N.s/m\\[0.5em]
$c_{4}$	&	102.625 N.s/m &	122.734 N.s/m\\[0.5em]
\bottomrule
\end{tabular}
\end{center}
\end{table}

\newpage
\section{Conclusion}
\label{sec:conclusion}

The current paper makes the following contributions in the area of design optimization of dynamic flexible MBSs. First, it presents the novel application of the DAVM for computing the sensitivities in design optimization of flexible dynamic multibody systems. In order to derive the adjoint equations, the discrete, instead of continuous, dynamic equations are used directly. This leads to a system of linear algebraic equations, which is in general less computationally expensive to solve compared to the system of differential algebraic equations produced by the CAVM. It also yields the exact values of gradients of discrete objective and constraint functions, thus ensuring the move along the actual descent or ascent directions in the optimization problems. Although the relevant adjoint equations are generated using a specific type of geometric variational integrators for the forward simulation, the proposed sensitivity analysis technique can be integrated, with minor modifications, with other types of numerical time-stepping methods as well.

Another contribution is demonstrating how to handle geometrical, as well as non-geometrical, design variables. From the design perspective, geometrical variables such as initial joint positions and linkage lengths are challenging yet significant part of the design embodiment process. The inclusion of the geometrical variables increases the design space dimensionality, and subsequently assists the designer in finding an optimum solution. A key consideration in dealing with geometrical design variables is the feasibility of the generated solutions. To prevent creating infeasible designs, a proper set of constraints should be determined for the optimization routine. For example, in the design of a rigid four-bar linkage where the crank is needed to make a full revolution, relevant dimensional constraints (according to Grashof's rules) must be devised to satisfy this requirement.

Lastly, this paper shows how rigid and flexible bodies (e.g., beams) can be effectively incorporated by utilizing ANCF and natural coordinates. The use of rotation-free formulations results in constant mass matrices for the bodies, which extensively simplifies the equations involved in both forward and backward simulations. In addition, it facilitates differentiating the equations with respect to geometrical and non-geometrical design variables, hence serving as an essential component for the contributions described above.

{It is worth noting that the proposed sensitivity analysis technique guarantees that the optimization routine finds an optimum solution and moves toward improving the defined performance criteria, as proved in Section {\ref{sec:problemdef}} and shown by the three numerical examples. It is well-known that if the initial solution is far away from the global optimum, it is possible that the optimization scheme falls into a local optimum. This is an intrinsic characteristic of all gradient-based optimizers and has been addressed extensively in the mathematics and engineering literature\mbox{\cite{NoceWrig06, gill1981practical}}. A possible remedy to this issue is to start the optimization procedure from different initial designs, thus increasing the chance of finding the global optimum.}

As a closing remark, the authors would like to emphasize that with the growing interests in generative design (i.e., design automation) in engineering applications, it is imperative to develop simulation and numeral techniques that can be effective in searching high-dimensional design spaces. The current work has highlighted some of the challenges in this endeavor and contributed an effective solution in the domain of dynamic multibody systems.

\begin{appendices}
\newpage
\section{Beam's mass and stiffness matrices}
\label{app:beammass}

Differentiating Equation \ref{eq:qs} with respect to time results in the global velocity vector below for a beam element

\begin{equation}
\dot{\vec{\mathbf{r}}}(x)=\mathbf{S}(x)\dot{\vec{\mathbf{q}}}(t)
\label{eq:dot}
\end{equation}

Using this and the kinetic energy equation, the element’s mass matrix can be obtained as

\begin{equation}\label{eq:kineticapp}
\begin{aligned}
T_{\text{element}} &= \frac{1}{2} \int\limits_V\left(\rho\left[\dot{\vec{\mathbf{r}}}^\mathrm{T}\dot{\vec{\mathbf{r}}}\right]\right) dV= \frac{1}{2} \dot{\vec{\mathbf{q}}}^\mathrm{T} \left[ \int\limits_V \left(\rho \left[\mathbf{S}^\mathrm{T}\mathbf{S} \right] \right) dV \right] \dot{\vec{\mathbf{q}}} = \frac{1}{2} \dot{\vec{\mathbf{q}}}^\mathrm{T} \mathbf{M}_{\text{element}} \dot{\vec{\mathbf{q}}} \\
\mathbf{M}_{\text{element}} &= \frac{\rho A l}{420} \begin{bmatrix}
156 \mathbf{I}_{3 \times 3} & 22l \mathbf{I}_{3 \times 3} & 54 \mathbf{I}_{3 \times 3} & -13l \mathbf{I}_{3 \times 3} \\[0.8em]
22l \mathbf{I}_{3 \times 3} & 4l^2 \mathbf{I}_{3 \times 3} & 13l \mathbf{I}_{3 \times 3} & -3l^2 \mathbf{I}_{3 \times 3} \\[0.8em]
54 \mathbf{I}_{3 \times 3} & 13l \mathbf{I}_{3 \times 3} & 156 \mathbf{I}_{3 \times 3} & -22l \mathbf{I}_{3 \times 3} \\[0.8em]
-13l \mathbf{I}_{3 \times 3} & -3l^2 \mathbf{I}_{3 \times 3} & -22l \mathbf{I}_{3 \times 3} & 4l^2 \mathbf{I}_{3 \times 3} \\
\end{bmatrix}
\end{aligned}
\end{equation}

$\mathbf{K}_{\text{longitudinal}}$ and $\mathbf{K}_{\text{transverse}}$ for the element are 

\begin{equation}\label{eq:linearFormapp}
\begin{aligned}
\mathbf{K}_{\text{longitudinal}} &= \frac{EA}{l} \varepsilon \begin{bmatrix}
\frac{6}{5} \mathbf{I}_{3 \times 3} & \frac{l}{10} \mathbf{I}_{3 \times 3} & -\frac{6}{5} \mathbf{I}_{3 \times 3} & \frac{l}{10} \mathbf{I}_{3 \times 3} \\[0.8em]
\frac{l}{10} \mathbf{I}_{3 \times 3} & \frac{2l^2}{15} \mathbf{I}_{3 \times 3} & -\frac{l}{10} \mathbf{I}_{3 \times 3} & -\frac{l^2}{30} \mathbf{I}_{3 \times 3} \\[0.8em]
-\frac{6}{5} \mathbf{I}_{3 \times 3} & -\frac{l}{10} \mathbf{I}_{3 \times 3} & \frac{6}{5} \mathbf{I}_{3 \times 3} & -\frac{l}{10} \mathbf{I}_{3 \times 3} \\[0.8em]
\frac{l}{10} \mathbf{I}_{3 \times 3} & -\frac{l^2}{30} \mathbf{I}_{3 \times 3} & -\frac{l}{10} \mathbf{I}_{3 \times 3} & \frac{2l^2}{30} \mathbf{I}_{3 \times 3}
\end{bmatrix} \\[0.8em]
\mathbf{K}_{\text{transverse}} &= \frac{EA}{l^3} \varepsilon \begin{bmatrix}
12 \mathbf{I}_{3 \times 3} & 6l \mathbf{I}_{3 \times 3} & -12 \mathbf{I}_{3 \times 3} & 6l \mathbf{I}_{3 \times 3} \\[0.8em]
6l \mathbf{I}_{3 \times 3} & 4l^2 \mathbf{I}_{3 \times 3} & -6l \mathbf{I}_{3 \times 3} & 2l^2 \mathbf{I}_{3 \times 3} \\[0.8em]
-12 \mathbf{I}_{3 \times 3} & -6l \mathbf{I}_{3 \times 3} & 12 \mathbf{I}_{3 \times 3} & -6l \mathbf{I}_{3 \times 3} \\[0.8em]
6l \mathbf{I}_{3 \times 3} & 2l^2 \mathbf{I}_{3 \times 3} & -6l \mathbf{I}_{3 \times 3} & 4l^2 \mathbf{I}_{3 \times 3} \\
\end{bmatrix}
\end{aligned}
\end{equation}

\section{Rigid body's mass matrix}
\label{app:rigid}

Using the kinetic energy equation, 
\begin{equation}\label{eq:kineticRigid}
T_{\text{rigid}} = \frac{1}{2} \int\limits_V\left(\rho\left[\dot{\vec{\mathbf{r}}}^\mathrm{T}\dot{\vec{\mathbf{r}}}\right]\right) dV=
\frac{1}{2} \dot{\vec{\mathbf{q}}}^\mathrm{T} \left[ \int\limits_V \left(\rho \left[\mathbf{S}^\mathrm{T}\mathbf{S} \right] \right) dV \right] =
\frac{1}{2} \dot{\vec{\mathbf{q}}}^\mathrm{T} \mathbf{M}_{\text{rigid}} \dot{\vec{\mathbf{q}}}
\end{equation}
and $\mathbf{S}$ in Equation \ref{eq:rp2}, the body's mass matrix becomes
\begin{equation}\label{eq:rigidmass}
\mathbf{M}_{\text{rigid}} = \int\limits_V \rho \begin{bmatrix}
\mathbf{I}_{3 \times 3} & \bar{x} \mathbf{I}_{3 \times 3} & \bar{y} \mathbf{I}_{3 \times 3} & \bar{z} \mathbf{I}_{3 \times 3} \\[0.8em]
\bar{x} \mathbf{I}_{3 \times 3} & \bar{x}^2 \mathbf{I}_{3 \times 3} & \bar{x}\bar{y} \mathbf{I}_{3 \times 3} & \bar{x}\bar{z} \mathbf{I}_{3 \times 3} \\[0.8em]
\bar{y} \mathbf{I}_{3 \times 3} & \bar{x}\bar{y} \mathbf{I}_{3 \times 3} & \bar{y}^2 \mathbf{I}_{3 \times 3} & \bar{y}\bar{z} \mathbf{I}_{3 \times 3} \\[0.8em]
\bar{z} \mathbf{I}_{3 \times 3} & \bar{x}\bar{z} \mathbf{I}_{3 \times 3} & \bar{y}\bar{z} \mathbf{I}_{3 \times 3} & \bar{z}^2 \mathbf{I}_{3 \times 3} \\
\end{bmatrix} dV
\end{equation}

If the origin of the body's local frame is located at its center of mass and the local frame's axes coincide with the body's principal axes, the mass matrix takes a diagonal form.
\section{Spring and damper's generalized force vectors}
\label{app:spring}

The potential energy of the spring depicted in Figure \ref{fig:case1} is defined by

\begin{equation}\label{eq:sprdampot}
U_{\mathrm{spring}}=\frac{1}{2}k(l_{t_{\mathrm{current}}}-l_0)^2
\end{equation}

in which $l_0$ and $l_{t_{\mathrm{current}}}$ are, respectively, the spring's lengths at zero and current time-step. Considering Equations \ref{eq:position} and \ref{eq:rp2}, Equation \ref{eq:sprdampot} can be rewritten as

\begin{equation}
\begin{aligned}
U_{\mathrm{spring}} & =\frac{1}{2}k\left( \sqrt{\left( \vec{\mathbf{r}}_{\mathrm{C}} -\vec{\mathbf{r}}_{\mathrm{D}} \right)^{\mathrm{T}} \left( \vec{\mathbf{r}}_{\mathrm{C}} -\vec{\mathbf{r}}_{\mathrm{D}} \right)}-l_0 \right)^2 \\
& = \frac{1}{2}k \left( \sqrt{ \vec{\bar{\mathbf{q}}}^{\mathrm{T}} \begin{bmatrix} \mathbf{S}_{\text{rigid}}^{\mathrm{T}} \Big|_{\mathrm{C}}  \mathbf{S}_{\text{rigid}} \Big|_{\mathrm{C}} & -\mathbf{S}_{\text{rigid}}^{\mathrm{T}} \Big|_{\mathrm{C}} \mathbf{S}_{\text{beam}} \Big|_{\mathrm{D}} \\[1.0em]
-\mathbf{S}_{\text{beam}}^{\mathrm{T}} \Big|_{\mathrm{D}}  \mathbf{S}_{\text{rigid}} \Big|_{\mathrm{C}} & \mathbf{S}_{\text{beam}}^{\mathrm{T}} \Big|_{\mathrm{D}} \mathbf{S}_{\text{beam}} \Big|_{\mathrm{D}}
\end{bmatrix} \vec{\bar{\mathbf{q}}} }-l_0 \right)^2
\end{aligned}
\end{equation}

where $\vec{\bar{\mathbf{q}}}=\begin{bmatrix} \vec{\mathbf{q}}_{\mathrm{rigid}}^{\mathrm{T}} & \vec{\mathbf{q}}_{\mathrm{beam}}^{\mathrm{T}} \end{bmatrix}^{\mathrm{T}}$. The \emph{generalized spring force} vector, thus, becomes

\begin{equation}
\frac{\partial U_{\mathrm{spring}}}{\partial \vec{\bar{\mathbf{q}}}}=k \frac{l_{t_{\mathrm{current}}}-l_0}{l_{t_{\mathrm{current}}}} \begin{bmatrix} \mathbf{S}_{\text{rigid}}^{\mathrm{T}} \Big|_{\mathrm{C}}  \mathbf{S}_{\text{rigid}} \Big|_{\mathrm{C}} & -\mathbf{S}_{\text{rigid}}^{\mathrm{T}} \Big|_{\mathrm{C}} \mathbf{S}_{\text{beam}} \Big|_{\mathrm{D}} \\[1.0em]
-\mathbf{S}_{\text{beam}}^{\mathrm{T}} \Big|_{\mathrm{D}}  \mathbf{S}_{\text{rigid}} \Big|_{\mathrm{C}} & \mathbf{S}_{\text{beam}}^{\mathrm{T}} \Big|_{\mathrm{D}} \mathbf{S}_{\text{beam}} \Big|_{\mathrm{D}}
\end{bmatrix} \begin{Bmatrix}
\vec{\mathbf{q}}_{\text{rigid}} \\[1.0em] \vec{\mathbf{q}}_{\text{beam}}
\end{Bmatrix}
\end{equation}

Following a similar procedure, the vector of \emph{generalized damper forces} is stated via

\begin{equation}\label{eq:dampereq}
\vec{\mathbf{f}}_{\text{damper}}= -c \begin{bmatrix} \mathbf{S}_{\text{rigid}}^{\mathrm{T}} \Big|_{\mathrm{C}}  \mathbf{S}_{\text{rigid}} \Big|_{\mathrm{C}} & -\mathbf{S}_{\text{rigid}}^{\mathrm{T}} \Big|_{\mathrm{C}} \mathbf{S}_{\text{beam}} \Big|_{\mathrm{D}} \\[1.0em]
-\mathbf{S}_{\text{beam}}^{\mathrm{T}} \Big|_{\mathrm{D}}  \mathbf{S}_{\text{rigid}} \Big|_{\mathrm{C}} & \mathbf{S}_{\text{beam}}^{\mathrm{T}} \Big|_{\mathrm{D}} \mathbf{S}_{\text{beam}} \Big|_{\mathrm{D}}
\end{bmatrix} \begin{Bmatrix}
\dot{\vec{\mathbf{q}}}_{\text{rigid}} \\[1.0em] \dot{\vec{\mathbf{q}}}_{\text{beam}}
\end{Bmatrix}
\end{equation}

\end{appendices}

\newpage
\bibliography{main}

\begin{thebibliography}{10}

\bibitem{Tromme2013}
E.~Tromme, O.~Br{\"{u}}ls, J.~Emonds-Alt, M.~Bruyneel, G.~Virlez, and
  P.~Duysinx, ``{Discussion on the optimization problem formulation of flexible
  components in multibody systems},'' {\em Structural and Multidisciplinary
  Optimization}, vol.~48, no.~6, pp.~1189--1206, 2013.

\bibitem{Greene1989}
W.~H. Greene and R.~T. Haftka, ``{Computational aspects of sensitivity
  calculations in transient structural analysis},'' {\em Computers {\&}
  Structures}, vol.~32, no.~2, pp.~433--443, 1989.

\bibitem{Maly1996}
T.~Maly and L.~R. Petzold, ``{Numerical methods and software for sensitivity
  analysis of differential-algebraic systems},'' {\em Applied Numerical
  Mathematics}, vol.~20, no.~1-2, pp.~57--79, 1996.

\bibitem{Mukherjee2008}
R.~M. Mukherjee, K.~D. Bhalerao, and K.~S. Anderson, ``{A divide-and-conquer
  direct differentiation approach for multibody system sensitivity analysis},''
  {\em Structural and Multidisciplinary Optimization}, vol.~35, no.~5,
  pp.~413--429, 2008.

\bibitem{Dias1997}
J.~M.~P. Dias and M.~S. Pereira, ``{Sensitivity analysis of rigid-flexible
  multibody systems},'' {\em Multibody System Dynamics}, vol.~1, no.~3,
  pp.~303--322, 1997.

\bibitem{Wang2005}
X.~Wang, E.~J. Haug, and W.~Pan, ``{Implicit numerical integration for design
  sensitivity analysis of rigid multibody systems},'' {\em Mechanics Based
  Design of Structures and Machines}, vol.~33, no.~1, pp.~1--30, 2005.

\bibitem{Haug1982}
E.~J. Haug and P.~E. Ehle, ``{Second-order design sensitivity analysis of
  mechanical system dynamics},'' {\em International Journal for Numerical
  Methods in Engineering}, vol.~18, no.~11, pp.~1699--1717, 1982.

\bibitem{Pi2012}
T.~Pi, Y.~Zhang, and L.~Chen, ``{First order sensitivity analysis of flexible
  multibody systems using absolute nodal coordinate formulation},'' {\em
  Multibody System Dynamics}, vol.~27, no.~2, pp.~153--171, 2012.

\bibitem{tromme2015structural}
E.~Tromme, D.~Tortorelli, O.~Br{\"u}ls, and P.~Duysinx, ``Structural
  optimization of multibody system components described using level set
  techniques,'' {\em Structural and Multidisciplinary Optimization}, vol.~52,
  no.~5, pp.~959--971, 2015.

\bibitem{Schaffer2005}
A.~S. Schaffer, ``{On the adjoint formulation of design sensitivity analysis of
  multibody dynamics},'' {\em Theses and Dissertations}, p.~93, 2005.

\bibitem{Ding2007}
J.-Y. Ding, Z.-K. Pan, and L.-Q. Chen, ``{Second order adjoint sensitivity
  analysis of multibody systems described by differential–algebraic
  equations},'' {\em Multibody System Dynamics}, vol.~18, no.~4, pp.~599--617,
  2007.

\bibitem{Zhu2015}
Y.~Zhu, D.~Dopico, C.~Sandu, and A.~Sandu, ``{Dynamic response optimization of
  complex multibody systems in a penalty formulation using adjoint
  sensitivity},'' {\em Journal of Computational and Nonlinear Dynamics},
  vol.~10, no.~3, p.~31009, 2015.

\bibitem{Cao2002}
Y.~Cao, S.~Li, and L.~Petzold, ``{Adjoint sensitivity analysis for
  differential-algebraic equations: algorithms and software},'' {\em Journal of
  Computational and Applied Mathematics}, vol.~149, no.~1, pp.~171--191, 2002.

\bibitem{Li2000}
S.~Li, L.~Petzold, and W.~Zhu, ``{Sensitivity analysis of
  differential–algebraic equations: a comparison of methods on a special
  problem},'' {\em Applied Numerical Mathematics}, vol.~32, no.~2,
  pp.~161--174, 2000.

\bibitem{ober2011discrete}
S.~Ober-Bl{\"o}baum, O.~Junge, and J.~E. Marsden, ``Discrete mechanics and
  optimal control: an analysis,'' {\em ESAIM: Control, Optimisation and
  Calculus of Variations}, vol.~17, no.~2, pp.~322--352, 2011.

\bibitem{Lauss2017a}
T.~Lau{\ss}, S.~Oberpeilsteiner, W.~Steiner, and K.~Nachbagauer, ``{The
  discrete adjoint method for parameter identification in multibody system
  dynamics},'' {\em Multibody System Dynamics}, pp.~1--14, 2017.

\bibitem{Leyendecker2008}
S.~Leyendecker, J.~E. Marsden, and M.~Ortiz, ``{Variational integrators for
  constrained dynamical systems},'' {\em ZAMM‐Journal of Applied Mathematics
  and Mechanics/Zeitschrift f{\"{u}}r Angewandte Mathematik und Mechanik},
  vol.~88, no.~9, pp.~677--708, 2008.

\bibitem{Marsden2001}
J.~E. Marsden and M.~West, ``{Discrete mechanics and variational
  integrators},'' {\em Acta Numerica}, vol.~10, pp.~357--514, 2001.

\bibitem{Kane2000}
C.~Kane, J.~E. Marsden, M.~Ortiz, and M.~West, ``{Variational integrators and
  the Newmark algorithm for conservative and dissipative mechanical systems},''
  {\em International Journal for Numerical Methods in Engineering}, vol.~49,
  no.~10, pp.~1295--1325, 2000.

\bibitem{Stern2006}
A.~Stern and M.~Desbrun, ``{Discrete geometric mechanics for variational time
  integrators},'' in {\em ACM SIGGRAPH 2006 Courses}, pp.~75--80, ACM, 2006.

\bibitem{Shabana2001}
A.~A. Shabana and R.~Y. Yakoub, ``Three dimensional absolute nodal coordinate
  formulation for beam elements: theory,'' {\em Journal of Mechanical Design},
  vol.~123, no.~4, pp.~606--613, 2001.

\bibitem{Uhlar2009}
S.~Uhlar and P.~Betsch, ``{A rotationless formulation of multibody dynamics:
  Modeling of screw joints and incorporation of control constraints},'' {\em
  Multibody System Dynamics}, vol.~22, no.~1, pp.~69--95, 2009.

\bibitem{Shabana1996}
A.~A. Shabana, ``{An absolute nodal coordinate formulation for the large
  rotation and deformation analysis of flexible bodies},'' {\em Technical
  Report, Department of Mechanical Engineering, University of Illinois at
  Chicago}, 1996.

\bibitem{Gerstmayr2006}
J.~Gerstmayr and A.~A. Shabana, ``{Analysis of thin beams and cables using the
  absolute nodal co-ordinate formulation},'' {\em Nonlinear Dynamics}, vol.~45,
  no.~1, pp.~109--130, 2006.

\bibitem{Garcia-Vallejo2003}
D.~Garcia-Vallejo, J.~L. Escalona, J.~Mayo, and J.~Dominguez, ``{Describing
  rigid-flexible multibody systems using absolute coordinates},'' {\em
  Nonlinear Dynamics}, vol.~34, no.~1, pp.~75--94, 2003.

\bibitem{Gerstmayr2013}
J.~Gerstmayr, H.~Sugiyama, and A.~Mikkola, ``{Review on the absolute nodal
  coordinate formulation for large deformation analysis of multibody
  systems},'' {\em Journal of Computational and Nonlinear Dynamics}, vol.~8,
  no.~3, p.~31016, 2013.

\bibitem{Yakoub2001}
R.~Y. Yakoub and A.~A. Shabana, ``Three dimensional absolute nodal coordinate
  formulation for beam elements: implementation and applications,'' {\em
  Journal of Mechanical Design}, vol.~123, no.~4, pp.~614--621, 2001.

\bibitem{Berzeri2000}
M.~Berzeri and A.~A. Shabana, ``{Development of simple models for the elastic
  forces in the absolute nodal co-ordinate formulation},'' {\em Journal of
  Sound and Vibration}, vol.~235, no.~4, pp.~539--565, 2000.

\bibitem{Takahashi1999}
Y.~Takahashi and N.~Shimizu, ``{Study on elastic forces of the absolute nodal
  coordinate formulation for deformable beams},'' in {\em ASME Proceedings of
  Design Engineering Technical Conference}, pp.~33--40, 1999.

\bibitem{DeJalon1987}
J.~G. de~Jalon, J.~Unda, A.~Avello, and J.~M. Jim{\'{e}}nez, ``{Dynamic
  analysis of three-dimensional mechanisms in “natural” coordinates},''
  {\em Journal of Mechanisms, Transmissions, and Automation in Design},
  vol.~109, no.~4, pp.~460--465, 1987.

\bibitem{DeJalon2012}
J.~G. {De Jalon} and E.~Bayo, {\em {Kinematic and dynamic simulation of
  multibody systems: the real-time challenge}}.
\newblock Springer Science {\&} Business Media, 2012.

\bibitem{Betsch2001}
P.~Betsch and P.~Steinmann, ``{Constrained integration of rigid body
  dynamics},'' {\em Computer Methods in Applied Mechanics and Engineering},
  vol.~191, no.~3, pp.~467--488, 2001.

\bibitem{DeJalon2007}
J.~G. de~Jal{\'{o}}n, ``{Twenty-five years of natural coordinates},'' {\em
  Multibody System Dynamics}, vol.~18, no.~1, pp.~15--33, 2007.

\bibitem{Garcia-Vallejo2008}
D.~Garcia-Vallejo, J.~Mayo, J.~L. Escalona, and J.~Dominguez,
  ``{Three-dimensional formulation of rigid-flexible multibody systems with
  flexible beam elements},'' {\em Multibody System Dynamics}, vol.~20, no.~1,
  pp.~1--28, 2008.

\bibitem{johnson2014nlopt}
S.~G. Johnson, ``The nlopt nonlinear-optimization package,'' 2014.

\bibitem{Sun2016}
J.~Sun, Q.~Tian, and H.~Hu, ``{Structural optimization of flexible components
  in a flexible multibody system modeled via ANCF},'' {\em Mechanism and
  Machine Theory}, vol.~104, pp.~59--80, 2016.

\bibitem{stolpe2018equivalent}
M.~Stolpe, A.~Verbart, and S.~Rojas-Labanda, ``The equivalent static loads
  method for structural optimization does not in general generate optimal
  designs,'' {\em Structural and Multidisciplinary Optimization}, vol.~58,
  no.~1, pp.~139--154, 2018.

\bibitem{NoceWrig06}
J.~Nocedal and S.~J. Wright, {\em Numerical Optimization}.
\newblock New York, NY, USA: Springer, second~ed., 2006.

\bibitem{gill1981practical}
P.~E. Gill, W.~Murray, and M.~H. Wright, ``Practical optimization,'' 1981.

\end{thebibliography}
\bibliographystyle{ieeetr}
\end{document}